\documentclass[11pt]{amsart}
\usepackage[]{amsmath, amsthm, amsfonts, verbatim, amssymb}
\usepackage[mathscr]{eucal}
\newcommand\B{{\bf B}}

\def\cD{{\mathscr D}}
\def\cE{{\mathscr E}}

\def\ocM{\overline{\mathscr M}}
\def\oM{\overline M}

\def\cH{\mathscr H}

\def\cI{{\mathcal I}}
\def\cR{{\mathscr R}}
\def\cG{{\mathscr G}}

\def\cT{{\mathscr T}}
\def\bR{{\mathbb R}}
\def\bZ{{\mathbb Z}}
\def\bQ{{\mathbb Q}}
\def\bC{{\mathbb C}}

\def\cT{{\mathscr T}}

\def\cP{{\mathscr P}}
\def\cC{{\mathscr C}}
\begin{document}
\newtheorem {theo}{Theorem}
\newtheorem {coro}{Corollary}
\newtheorem {lemm}{Lemma}
\newtheorem {rem}{Remark}
\newtheorem {defi}{Definition}
\newtheorem {ques}{Question}
\newtheorem {prop}{Proposition}
\def\spb{\smallpagebreak}
\def\mpb{\vskip 0.5truecm}
\def\bpb{\vskip 1truecm}
\def\wtM{\widetilde M}
\def\tM{\widetilde M}
\def\wtN{\widetilde N}
\def\tN{\widetilde N}
\def\tC{\widetilde C}
\def\tX{\widetilde X}
\def\tY{\widetilde Y}
\def\ti{\widetilde \iota}
\def\bs{\bigskip}
\def\ms{\medskip}
\def\ni{\noindent}
\def\td{\nabla}
\def\pd{\partial}
\def\b{\bullet}
\def\hol{$\text{hol}\,$}
\def\Log{\mbox{Log}}
\def\bfQ{{\bf Q}}
\def\Todd{\mbox{Todd}}
\def\bP{{\bf P}}
\def\dxi{d x^i}
\def\dxj{d x^j}
\def\dyi{d y^i}
\def\dyj{d y^j}
\def\dzi{d z^I}
\def\dzj{d z^J}
\def\ozi{d{\overline z}^I}
\def\ozj{d{\overline z}^J}
\def\oz1{d{\overline z}^1}
\def\oz2{d{\overline z}^2}
\def\oz3{d{\overline z}^3}
\def\sI{\sqrt{-1}}
\def\hol{$\text{hol}\,$}
\def\ok{\overline k}
\def\ol{\overline l}
\def\oJ{\overline J}
\def\oT{\overline T}
\def\oS{\overline S}
\def\oV{\overline V}
\def\oW{\overline W}
\def\oI{\overline I}
\def\oK{\overline K}
\def\oL{\overline L}
\def\oj{\overline j}
\def\oi{\overline i}
\def\ok{\overline k}
\def\oz{\overline z}
\def\om{\overline mu}
\def\on{\overline nu}
\def\oa{\overline \alpha}
\def\ob{\overline \beta}
\def\of{\overline f}
\def\og{\overline \gamma}
\def\ogamma{\overline \gamma}
\def\odelta{\overline \delta}
\def\otheta{\overline \theta}
\def\ophi{\overline \phi}
\def\opd{\overline \partial}
\def\oA{\overline A} 
\def\oB{\overline B}
\def\oC{\overline C}
\def\op{\overline D}
\def\oIq1{\oI_1\cdots\oI_{q-1}}
\def\oIq2{\oI_1\cdots\oI_{q-2}}
\def\op{\overline \partial}
\def\ua{{\underline {a}}}
\def\us{{\underline {\sigma}}}
\def\tor{{\mbox{tor}}}
\def\vol{{\mbox{vol}}}
\def\rank{{\mbox{rank}}}
\def\bp{{\bf p}}
\def\bk{{\bf k}}
\def\a{{\alpha}}
\def\tchi{\widetilde{\chi}}

\vskip5mm

\title[arithmetic fake projective spaces and Grassmannians] 
{Arithmetic fake projective spaces and arithmetic fake Grassmannians}
 \maketitle
{\centerline{\sc Gopal Prasad and Sai-Kee Yeung}}
\vskip5mm

\centerline{\it Dedicated to Robert P.\:Langlands on his 70th birthday}
\vskip7mm

\begin{center}
{\bf 1. Introduction}  
\end{center}
\vskip4mm

Let $r$ be a positive integer, and let $Z$ be a compact K\"ahler manifold of dimension $r$ whose Betti numbers are same as that of  $\bP_{\bC}^{r}$ but which is not isomorphic to $\bP_{\bC}^{r}$. Then for $r=2$, since $c_1^2 = 9= 3c_2$, it follows from S.-T.\,Yau's results on Calabi's conjecture that $Z$ is uniformized by the open unit ball $\B^2$ in $\bC^2$, i.e., it is the quotient of $\B^2$ by a cocompact {\it torsion-free} discrete subgroup $\Pi$ of  the automorphism group $\mathrm{PU}(2,1)$ of $\B^2$. In this case (i.e., if $r=2$), it was proved further by Bruno Klingler [Kl], and the second author in [Y], that $\Pi$ is an arithmetic subgroup of $\mathrm{PU}(2,1)$. In this paper the quotient  of the open unit ball $\B^{r}$ in $\bC^{r}$ by a cocompact {\it torsion-free arithmetic} subgroup of the group $\mathrm{PU}(r,1)$ of automorphisms of $\B^{r}$ will be called an {\it arithmetic fake $\bP^{r}_{\bC}$} if it has the same Betti numbers as $\bP^{r}_{\bC}$. 

\vskip1mm

We observe that $\B^{r}$ is the symmetric space of ${\rm PU}(r,1)$, and $\bP_{\bC}^{r}$ is its compact dual. Now given a connected semi-simple real algebraic group $\overline{G}$ with trivial center, let $X$ be the symmetric space of $\overline{G}(\bR)$ ($X$ is the space of maximal compact subgroups of $G(\bR)$) and $X_u$ be the compact dual of $X$. We shall say that the quotient  $X/\Pi$ of $X$ by a cocompact {\it torsion-free arithmetic} subgroup $\Pi$ of $\overline{G}(\bR)$ is an {\it arithmetic fake $X_u$} if its Betti numbers are same as that of $X_u$; $X/\Pi$ is an {\it irreducible} arithmetic fake $X_u$ if, further, $\Pi$ is irreducible (i.e., no subgroup of $\Pi$ of finite index is a direct product of two infinite normal subgroups).  For example, the Grassmann space 
${\bf Gr}_{m,n}$ of $m$-dimensional vector subspaces of ${\bC}^n$ is the compact dual of the symmetric space of the group ${\rm PU}(n-m,m)$, and so the quotient of the symmetric space of ${\rm PU} (n-m,m)$ by a cocompact torsion-free arithmetic subgroup of ${\rm PU}(n-m,m)$, whose Betti numbers are same as that of ${\bf Gr}_{m,n}$, is an {\it irreducible arithmetic fake} ${\bf Gr}_{m.n}$.

\vskip1mm

Let $\overline{G}$, $X$ and $X_u$ be as above, and let $\Pi$ be a cocompact torsion-free arithmetic  subgroup of $\overline{G}(\bR)$. Let $Z = X/\Pi$.  If $Z$ is an (arithmetic) fake $X_u$, then the Euler-Poincar\'e characteristic $\chi(Z)$ of $Z=X/\Pi$, and so the Euler-Poincar\'e characteristic $\chi(\Pi)$ of 
$\Pi$, equals that of $X_u$. Using the results of [BP], we can easily conclude that there are only finitely many irreducible arithmetic fake $X_u$ with $\chi(X_u)\ne 0$ (in this finiteness assertion, $\overline{G}$ is allowed to vary). It is of interest to determine them all. If $\Pi$ is contained in the identity component of $\overline{G}(\bR)$, then there is a natural embedding of $H^*(X_u,\bC)$ in $H^*(Z,\bC)$; see, for example, [B], 3.1 and 10.2. Thus $Z=X/\Pi$ is an arithmetic fake $X_u$ if and only if the natural homomorphism $H^*(X_u,\bC)\to H^*(Z,\bC)$ is an isomorphism.  This latter property makes such arithmetic fake $X_u$ very interesting geometrically as well as  for the theory of automorphic forms.
\vskip 1mm

If the symmetric space $X$ is hermitian, then $Z$ is a smooth complex projective algebraic variety. Hence, if $X$ is hermitian, an arithmetic fake $X_u$ is a smooth complex projective algebraic variety. 

\vskip1mm

Let $n$ be an integer $>1$. The Euler-Poincar\'e characteristic of ${\bP}_{\bC}^{n-1}$, and so also of any arithmetic fake ${\bP}_{\bC}^{n-1}$, is $n$.  It is an immediate consequence of the Hirzebruch proportionality principle, see [Se1], Proposition 23, that the orbifold Euler-Poincar\'e characteristic (i.e., the Euler-Poincar\'e characteristic in the sense of C.T.C.\,Wall, see [Se1], \S 1.8) of any cocompact discrete subgroup of ${\rm PU}(n-1,1)$, for $n$ even, is negative. This implies that  if there exists an arithmetic fake ${\bP}_{\bC}^{n-1}$, then $n$ is necessarily {\it odd}.   The purpose of this paper is to determine all irreducible cocompact {\it torsion-free arithmetic} subgroups $\overline\Gamma$ of a product $\overline{\mathcal{G}}$ of $r$ groups of the form ${\rm PU}(n-m,m)$, with $n>3$ odd, and $0<m<n$,  whose Euler-Poincar\'e characteristic $\chi({\overline\Gamma})$ is equal to the Euler-Poincar\'e characteristic $\chi(X_u)$ of the compact dual $X_u$ of the!
  symmetric space $X$ of $\overline{\mathcal{G}}$.  (Note that $\chi(X_u)>0$.)
\vskip1mm
 
Let $\overline\Gamma$ be an irreducible cocompact torsion-free arithmetic subgroup of $\overline{\mathcal{G}}$ with $\chi({\overline\Gamma})= \chi(X_u )$. Let $\mathcal{G}$ be the connected semi-simple Lie group obtained by replacing each of the 
$r$ factors  ${\rm PU}(n-m,m)$ of ${\overline{\mathcal{G}}}$ by ${\rm SU}(n-m,m)$. As the kernel of the 
natural surjective homomorphism ${\mathcal{G}}\rightarrow {\overline{\mathcal{G}}}$ is a group of order $n^r$, if $\widetilde\Gamma$ is the full inverse image of $\overline\Gamma$ in $\mathcal{G}$, then $\widetilde\Gamma$ is an 
arithmetic subgroup whose orbifold Euler-Poincar\'e characteristic is $\chi(X_u)/n^r$. Therefore, the 
orbifold Euler-Poincar\'e characteristic of any arithmetic subgroup of $\mathcal{G}$, which contains $\widetilde\Gamma$, is a submultiple of $\chi(X_u)/n^r$. Assume, if possible, that $\mathcal{G}$ contains an irreducible  maximal arithmetic subgroup $\Gamma$ whose orbifold  
Euler-Poincar\'e characteristic $\chi(\Gamma)$ is a submultiple of $\chi(X_u)/n^r$.       
As $\Gamma$ is an irreducible maximal arithmetic subgroup of $\mathcal{G}$, there exist 
a totally real number field $k$ of degree at least $r$ over $\bQ$, an absolutely simple simply connected group $G$, of type $^2A_{n-1}$, defined 
over $k$, $r$ real places of $k$, say $v_j$, $j=1,\ldots,\,r$, such that   $\mathcal{G}\cong \prod_{j=1}^r G(k_{v_j})$,  and for every other real place $v$ of $k$, $G(k_v)$ is 
isomorphic to the compact real Lie group ${\rm SU}(n)$, 
and a ``principal'' arithmetic subgroup $\Lambda$ of $G(k)$ such that $\Gamma$ is the normalizer of 
$\Lambda$ in $\mathcal{G}$ (we identify $\mathcal{G}$ with $\prod_{j=1}^r G(k_{v_j})$ and use this identification to realize $G(k)$ as a subgroup of $\mathcal{G}$), see Proposition 1.4(iv) of [BP]. 
\vskip1mm

     From the description of absolutely simple groups of type $^2A_{n-1}$ 
(see, for example, [T1]), we know that there exists a quadratic extension $\ell$ of $k$, 
a division algebra $\cD$ with 
center $\ell$ and of degree $s =\sqrt{[\cD:\ell]}$, $s|n$, $\cD$ given 
with an involution $\sigma$ of the second kind such that $k=\{ x\in \ell \,
|\, x =\sigma(x) \}$, and a nondegenerate hermitian form $h$ on 
$\cD^{n/s}$ defined in terms of the involution $\sigma$ so that $G$ is 
the special unitary group ${\rm SU}(h)$ of $h$. It is obvious that 
$\ell$ is totally complex.

\vskip1mm

In terms 
of the normalized Haar-measure $\mu$ on $\mathcal{G}= \prod_{j=1}^r G(k_{v_j})$ used in [P]
and [BP], and to be used in this paper, $\chi(\Gamma)=\chi(X_u)\mu(\mathcal{G}/\Gamma)$ (see [BP], 4.2). Thus the condition that $\chi(\Gamma)$ is a submultiple of  $\chi(X_u)/n^r$ is equivalent to 
the condition that $\mu(\mathcal{G}/\Gamma)$ is a submultiple of ${1/n^r}.$ We shall prove that if $n>7$, there does not exist an arithmetic subgroup whose covolume is $\leqslant 1/n^r$, and if $n=5$ or $7$, there does not exist an arithmetic subgroup whose covolume  
is a submultiple of $1/n^r$. 
\vskip1mm

The main result (Theorem 2) of this paper implies that {\it arithmetic fake 
$\bP_{\bC}^{n-1}$ can exist only if  $n=3$ or $5$}, and {\it an arithmetic 
fake ${\bf Gr}_{m,n}$ exists, with $n>3$ odd, only if $n =5$}. The first 
fake projective plane was constructed by David Mumford in [M] using $p$-adic 
uniformization.  We have constructed twenty three distinct (finite) 
classes of arithmetic fake projective planes, and it has been proved that there 
can exist at most three more (see [PY], and the addendum to [PY]). In \S 5 of this paper we construct 
four distinct $4$-dimensional arithmetic fake projective spaces and four 
distinct fake ${\bf Gr}_{2,5}$. In \S 6, certain results and 
computations of [PY] have been used to exhibit five irreducible arithmetic fake $\bP^2_{\bC}\times \bP^2_{\bC}$. 
All these are connected smooth (complex projective) Shimura varieties, 
and these are the first examples of fake $\bP_{\bC}^4$, fake Grassmannians, 
and irreducible fake $\bP^2_{\bC}\times\bP^2_{\bC}$.

\vskip6mm
\begin{center}{\bf 2. Preliminaries}
\end{center}
\vskip1mm

A comprehensive survey of the basic definitions and the main results of the 
Bruhat--Tits theory of reductive groups over nonarchimedean local fields is given in [T2].
\vskip2mm

\ni{\bf 2.1.} Throughout this paper we will use the notations introduced in 
the introduction. $n$ will always be an odd integer $>3$, $k$ a totally real number field of degree $d$,   $\ell$ a totally complex quadratic extension of $k$, and $G = {\rm SU}(h)$, where $h$ is as in the introduction. $G$ is an absolutely simple simply connected $k$-group of type $^{2}A_{n-1}$. All unexplained notations are as in [BP] and [P]. Thus for a number field $K$, $D_K$ will denote the absolute value of its 
discriminant, $h_K$ its class number, i.e., the order of its class group 
$Cl (K)$. We will denote by $h_{K,n}$ the order of the subgroup of $Cl(K)$ consisting 
of the elements of order dividing $n$. Then $h_{K,n}
\leqslant h_K$. We will denote by $U_K$ the multiplicative-group of units of $K$, and
by $K_n$ the subgroup of $K^{\times}$ consisting of the elements $x$
such that for every normalized valuation $v$ of $K$, $v(x)\in n\bZ$.
\vskip1mm

 $V_f$ (resp.\,\,$V_{\infty}$) will denote the set of nonarchimedean
(resp.\,\,archimedean) places of $k$. As $k$ admits at least $r$ distinct real places, see the introduction, $d\geqslant r$.  For $v\in V_f$, $q_v$ will
denote 
the cardinality of 
the residue field ${\mathfrak f}_v$ of $k_v$. 
\vskip1mm

For a parahoric subgroup $P_v$ 
of $G(k_v)$,  we define $e(P_v)$ and $e'(P_v)$ by the following formulae (cf.\,Theorem 3.7 of [P]):
\begin{equation}
e(P_v)=\frac{q_v^{(\dim\oM_v+\dim\ocM_v)/2}}{\#\oM_v({\mathfrak f}_v)}.
\end{equation} 
\begin{equation}
e'(P_v) =e(P_v)\cdot \frac{\#\ocM_v({\mathfrak f}_v)}{q_v^{\dim\ocM_v}}=q_v^{(\dim\oM_v -\dim\ocM_v)/2}\cdot\frac{\#\ocM_v({\mathfrak f}_v)}{\#\oM_v({\mathfrak f}_v)}.
\end{equation}
If $v$ splits in $\ell$, then  
$$e'(P_v) = e(P_v)\prod_{j=1}^{n-1}(1-\frac{1}{q^{j+1}_v}),$$ 
and if $v$ does not split in $\ell$, then   
$$e'(P_v) = e(P_v)\prod_{j=1}^{(n-1)/2}(1-\frac{1}{q_v^{2j}})(1+\frac{1}
{q_v^{2j+1}}),$$ or $$e'(P_v) = e(P_v)\prod_{j=1}^{(n-1)/2}
(1-\frac{1}{q_v^{2j}})$$
according as $v$ {\it does not} or {\it does} ramify in $\ell$. It is obvious that $e'(P_v) < e(P_v)$, 
and it follows from Proposition 2.10 (iii) of [P] that for any parahoric subgroup $P'_v$ contained in $P_v$, $e'(P'_v)=[P_v:P'_v]e'(P_v)$. 

\vskip2mm

\ni{\bf 2.2.} We note that if $P_v$ is a {\it hyperspecial} parahoric
subgroup 
of $G(k_v)$, then the ${\mathfrak f}_v$-group ${\oM}_v$, which in 
this case is just the ``reduction mod $\mathfrak{p}$" of $P_v$,  is 
either ${\rm SL}_n$ or 
${\rm SU}_n$ according as $v$ does or does not split in $\ell$, and ${\ocM}_v ={\oM}_v$ . If $v$
ramifies in $\ell$, then $G$ is quasi-split over $k_v$, and if $P_v$
is {\it special}, then ${\oM}_v$ is isogenous
to either ${\rm SO}_{n}$ or ${\rm Sp}_{n-1}$, and so is $\ocM_v$. Hence, $e'(P_v) = 1$ if either $P_v$ is
hyperspecial, or $v$ ramifies in $\ell$ and $P_v$ is special. 
\vskip1mm

\ni{\bf 2.3.} (i) Let $v$ be a nonarchimedean place of $k$ which splits in $\ell$
and $G$ splits at $v$. Then $G$ is isomorphic to ${\rm SL}_n$ over $k_v$, and 
$\ocM_v$ is ${\mathfrak f}_v$-isomorphic
to ${\rm SL}_n$. It can be
seen by a 
direct computation that for any nonhyperspecial parahoric subgroup
$P_v$ 
of $G(k_v)$,
$e'(P_v)$ is an integer greater than $n$.    
\vskip1mm

{(ii)} Let $v$ be a nonarchimedean place of $k$ which splits in $\ell$ 
but $G$ does not split at $v$. Then $k_v\otimes_k \cD = ({k_v\otimes_k \ell})\otimes_{\ell}\cD = M_{n/d_v}({\mathfrak D}_v)\oplus \sigma(M_{n/d_v}({\mathfrak D}_v))$, where ${\mathfrak D}_v$ is a 
division algebra with center  $k_v$, of degree $d_v >1$,
$d_{v}|n$. Hence, $G$ is $k_v$-isomorphic to ${\rm
  SL}_{n/d_{v},{\mathfrak D}_v}$. Let $P_v$ be a
maximal parahoric subgroup of $G(k_v)$. Then $\ocM_v$ is 
${\mathfrak f}_v$-isomorphic to ${\rm SL}_n$, and  ${\overline M}_v$ 
is isogenous to the product of the norm-$1$ torus $R^{(1)}_{F_v/{\mathfrak f}_v}({\rm GL}_1)$ and the semi-simple group $R_{F_v/{\mathfrak f}_v}
({\rm SL}_{n/d_v})$, where $F_v$ is the field extension of ${\mathfrak f}_v$ of degree $d_v$. So $$\#{\overline M}_v({\mathfrak f}_v )= \frac{q_v^{n^2 /d_{v}}}{q_v -1}\prod_{j =1}^{n/d_v}\big( 1-\frac{1}{q_v^{jd_v}}
\big),$$
and hence, $$e'(P_v) = q_v^{n^2(d_v-1)/2d_v}
\prod_{j=1}^{n/d_v}\big(1-\frac{1}{q_v^{jd_v}}\big)^{-1}\prod_{j=1}^{n}
\big( 1-\frac{1}{q_v^{j}}\big)$$ 
$$\ \ \ \ \ \ \ \ \ \ =
 \frac{\prod_{j=1}^{n}(q_v^j -1)}{\prod_{j=1}^{n/d_v}(q_v^{jd_v}-1)}> q_v^{(n^2-2n)(d_v-1)/2d_v}>n.$$

The above computation shows also that $e'(P_v)$ is an integer.
 \vskip1mm
 
 {(iii)} Now let $v$ be a nonarchimedean place of $k$ which does not split in $\ell$. Then being a group of type $^2A_{n-1}$, with $n$ odd, $G$ is quasi-split over $k_v$. It is not difficult to see, using (2), the fact that the order of a subgroup of a finite group divides the order of the latter, an obvious analogue for connected reductive algebraic groups defined over a finite field of the results of Borel and de Siebenthal [BdS] on subgroups of maximal rank of a compact connected Lie group, and the fact that for a finite field $\mathfrak f$, the groups of $\mathfrak f$-rational points of connected absolutely simple $\mathfrak f$-groups of type $B_m$ and $C_m$, for an arbitrary $m$,  have equal order, that  $e'(P_v)$ is an integer.
\vskip1mm
 
>From (i), (ii) and (iii) we gather that {\it for all $v\in V_f$, $e'(P_v)$ is an integer.}

\vskip2mm

 \ni{\bf 2.4.} Let $\Gamma$ be a maximal arithmetic subgroup of
$\mathcal{G}=\prod_{j=1}^rG(k_{v_j})$ such that $n^r\mu (\mathcal{G}/\Gamma)\leqslant 1$, see the introduction. Let
$\Lambda = \Gamma\cap G(k)$. Then $\Gamma$ is the normalizer of $\Lambda$ in $\mathcal{G}$, 
and $\Lambda$ is a {\it principal} arithmetic
subgroup (see [BP], 
Proposition 1.4(iv)), i.e., for every nonarchimedean place $v$ of $k$, 
the closure $P_v$ of $\Lambda$ in $G(k_v)$ is a parahoric
subgroup, and $\Lambda = G(k)\cap \prod_{v\in V_f}P_v$. Let $\cT$ be
the set of all nonarchimedean $v$ which splits in $\ell$ and
$P_v$ is not a hyperspecial parahoric subgroup of $G(k_v)$. Let $\cT^{\prime}$
be the set of all nonarchimedean $v$ which does not split in
$\ell$, and either $P_v$ is not a hyperspecial parahoric subgroup of
$G(k_v)$ but a hyperspecial parahoric exists (which is the case if, and only if, 
$v$ is unramified over $\ell$), or $v$ is ramified in $\ell$ and $P_v$
is not a special parahoric subgroup.  
\vskip1mm 

\ni{\bf 2.5.} Let $\mu_n$ be the kernel of the endomorphism $x\mapsto
x^n$ 
of ${\rm  GL}_1$. Then the center $C$ of $G$ is $k$-isomorphic to the kernel of the norm map 
$N_{\ell/k}$ from the algebraic group $R_{\ell/k}(\mu_n)$, obtained
from $\mu_n$  by Weil's restriction of scalars, to $\mu_n$. 

As $n$ is odd,  the norm map $N_{\ell/k}:\, \mu_n(\ell)\rightarrow
\mu_n(k)$ is onto, $\mu_n(k)/N_{\ell/k}(\mu_n(\ell))$ is trivial, and
hence, the Galois cohomology group $H^1(k,C)$ is isomorphic to the
kernel of the homomorphism $\ell^{\times}/{\ell^{\times}}^n\rightarrow
k^{\times}/{k^{\times}}^n$ induced by the norm map. We shall denote
this kernel by $(\ell^{\times}/{\ell^{\times}}^n)_{\b}$ in the sequel.

By Dirichlet's unit theorem, $U_k\cong \{\pm 1\}\times 
{\bZ}^{d-1}$, and $U_{\ell}\cong \mu(\ell)\times {\bZ}^{d-1}$, where
$\mu(\ell)$ is the finite cyclic group of roots of unity in $\ell$. Hence, 
$U_k/U_k^n\cong (\bZ/n\bZ)^{d-1}$, and $U_{\ell}/U_{\ell}^n\cong
\mu(\ell)_n\times (\bZ/n\bZ)^{d-1}$, where $\mu(\ell)_n$ is the group
of $n$-th roots of unity in $\ell$. Now we observe that
$N_{\ell/k}(U_{\ell})\supset N_{\ell/k}(U_k) = U_k^2$, which implies
that, as $n$ is odd, the homomorphism  $U_{\ell}/U_{\ell}^n\rightarrow 
U_k/U_k^n$, induced by the norm map, is onto. Therefore, the order of
the kernel $(U_{\ell}/U_{\ell}^n)_{\b}$ of this homomorphism equals $\#
\mu(\ell)_n$. 

\vskip1mm
 
The short exact sequence $(4)$ in the proof of Proposition 0.12 of  [BP]
gives us the following exact sequence: $$1\rightarrow
(U_{\ell}/U_{\ell}^n)_{\b}\rightarrow
(\ell_n/{\ell^{\times}}^n)_{\b}\rightarrow (\cP\cap \cI^n)/\cP^n,$$
where $(\ell_n/{\ell^{\times}}^n)_{\b} = (\ell_n/{\ell^{\times}}^n)\cap
(\ell^{\times}/{\ell^{\times}}^n)_{\b}$, $\cP$ is the group of all
fractional 
principal ideals of $\ell$, and $\cI$ the group of all fractional ideals (we use multiplicative
notation for the group operation in both $\cI$ and $\cP$). Since the
order of the last group of the above exact sequence is $h_{\ell,n}$,
see $(5)$ in the proof of Proposition 0.12 of  [BP], we conclude that 
$$\#(\ell_n/{\ell^{\times}}^n)_{\b} \leqslant \#\mu(\ell)_n\cdot h_{\ell,n}.$$

Now we note that the order of the first term of the short exact sequence 
of Proposition 2.9 of [BP], for $G' =G$ and $S =V_{\infty}$, 
is $n^r /\#\mu(\ell)_n$.

Using the above observations, together with Proposition 2.9 and Lemma
5.4 
of [BP], and a
close 
look at the 
arguments in 5.3 and 5.5 of [BP] for $S=V_{\infty}$ and $G$ as above, we can derive the following upper bound: 
\begin{equation}
[\Gamma : \Lambda] \leqslant n^{r+\#\cT}h_{\ell, n}.  
\end{equation}

\vskip1mm

\ni From this we obtain
\begin{equation} 1 \geqslant n^r\mu(\mathcal{G}/\Gamma)\geqslant \frac{\mu(\mathcal{G}/\Lambda)}{n^{\#\cT}h_{\ell,n}}.
\end{equation}

\vskip1mm

\ni{\bf 2.6.} Now we will use the volume formula of [P] to write down
the precise value of $\mu(\mathcal{G}/\Lambda)$. As the Tamagawa number 
$\tau_k(G)$ of $G$ equals $1,$ Theorem 3.7 of [P] (recalled in 3.7 of [BP]), for $S=V_{\infty}$,  gives us for $n$ odd, 
\begin{equation}
\mu(\mathcal{G}/\Lambda )= D_k^{(n^2-1)/2}(D_\ell/ D_k^2)^{(n-1)(n+2)/4}\Big(\prod_{j=1}^{n-1}\frac{j!}{(2\pi)^{j+1}}\Big)^{d}\cE, 
\end{equation}
where
$\cE=\prod_{v\in V_f}e(P_v),$ with $e(P_v)$ as in 2.1.

\vskip2mm

\ni{\bf 2.7.} Let $\zeta_k$ be the Dedekind zeta-function of $k$, and $L_{\ell | k}$ be the Hecke $L$-function associated to the quadratic Dirichlet character of $\ell/k$. Then
$$\zeta_k(j) =\prod_{v\in V_f}(1-\frac1{q_v^j})^{-1};$$  
$$L_{\ell|k}(j) ={\prod} '(1-\frac1{q_v^j})^{-1}{\prod} ''
(1+\frac1{q_v^j})^{-1},$$
where $\prod'$ is the product over the nonarchimedean places $v$ of $k$ which 
split in $\ell$, and $\prod''$
is the product over all the other nonarchimedean places $v$ which do not 
ramify in $\ell$. Hence the Euler product $\cE$ appearing in (5) 
can be rewritten as
\begin{equation}
\cE=\prod_{v\in V_f} e^{\prime}(P_v) \prod_{j=1}^{(n-1)/2}\big(
\zeta_k(2j)L_{\ell|k}(2j+1)\big).$$ 
Since $e'(P_v) = 1$, if $v\notin \cT\cup \cT'$ (2.2), and $e'(P_v)>n$ if $v\in \cT$ (2.1 and 2.3),  
$$\cE = \prod_{v\in\cT\cup\cT'} e'(P_v) \prod_{j=1}^{(n-1)/2}\big(\zeta_k(2j) 
L_{\ell|k}(2j+1)\big)
\geqslant n^{\#\cT}\prod_{j=1}^{(n-1)/2}\big(\zeta_k(2j) 
L_{\ell|k}(2j+1)\big).
\end{equation}

\ni{\bf 2.8.} Using the functional equations $$\zeta_k(2j) = D_k^{\frac{1}{2}-
2j}\big(\frac{(-1)^{j}2^{2j-1}\pi^{2j}}{(2j-1)!}\big)^d \zeta_k(1-2j),$$
and $$L_{\ell|k}(2j+1) = \big(\frac{D_k}{D_{\ell}}\big)^{2j+\frac{1}{2}}\big(\frac{(-1)^j 2^{2j}\pi^{2j+1}}{(2j)!}\big)^d L_{\ell|k}(-2j),$$
we find that 
\begin{equation}
\cR :=D_k^{(n^2-1)/2}(D_\ell/ D_k^2)^{(n-1)(n+2)/4}\Big(\prod_{j=1}^{n-1}\frac{j!}{(2\pi)^{j+1}}\Big)^{d}\prod^{(n-1)/2}_{j=1}\big(     
\zeta_k(2j)L_{\ell|k}(2j+1)\big)
\end{equation}

$$\ \ \ \ =2^{-(n-1)d}\zeta_k(-1)L_{\ell|k}(-2)\zeta_k(-3)L_{\ell|k}(-4)\cdots \zeta_k(2-n) 
L_{\ell|k}(1-n).$$
Equations (5), (6) and (7) imply that  $$\mu(\mathcal{G}/\Lambda)= \cR\prod_{\cT \cup \cT'}e'(P_v).$$ 
As $e'(P_v)$ is an integer for all $v$ (see 2.3), we conclude that $\mu(\mathcal{G}/\Lambda)$ is an integral multiple of $\cR$.  
\vskip2mm

\ni{\bf 2.9.}  As $\chi(\Lambda)= \chi(X_u)\mu(\mathcal{G}/{\Lambda})$ ([BP], 4.2), we have the following\begin{equation}
\chi(\Gamma)= \frac{\chi(\Lambda)}{[\Gamma :\Lambda ]}= \frac{\chi(X_u)\mu (\mathcal{G}/{\Lambda})}
{[\Gamma : \Lambda]}.
\end{equation}
Proposition 2.9 of [BP] applied to $G' = G$ and $\Gamma' = \Gamma$ implies that any prime divisor of 
the integer $[\Gamma :\Lambda]$ divides $n$. Now since $\mu(\mathcal{G}/\Lambda)$ is an integral multiple of $\cR$ (the latter as in (7)), we 
conclude from (8) that  if $\chi(\Gamma)$ is a submultiple of $\chi(X_u)$, then any prime which divides the numerator of the rational number 
$\cR$ is a divisor of $n$. We state this as the following proposition.
\begin{prop}
If the orbifold Euler-Poincar\'e characteristic of $\Gamma$ is a submultiple of $\chi(X_u)$, then any prime divisor of the  
numerator of the rational number $\cR$ divides $n$.
\end{prop}

\vskip2mm

\ni{\bf 2.10.} We know (cf.\,\,[P], Proposition 2.10(iv), and 2.3 above) that 
\begin{equation}
 {\rm for \ all} \ v\in V_f, \ \ e(P_v)>1, \ {\rm and \ for\ all}\ v\in \cT,\ \  e(P_v)> e'(P_v)>n.
\end{equation}

\ni Now combining (4), (5) and (9), we obtain 
\begin{equation}
 1\geqslant n^r\mu(\mathcal{G}/\Gamma) > \frac{{D_\ell}^{(n-1)(n+2)/4}}{D_k^{(n-1)/2}h_{\ell,n}}\Big(\prod_{j=1}^{n-1}\frac{j!}{(2\pi)^{j+1}}\Big)^{d}.
\end{equation}
\vskip2mm

\ni It follows 
from Brauer-Siegel Theorem that for all real $s>1$, 
\begin{equation}
 h_\ell R_\ell\leqslant w_\ell s(s-1)\Gamma(s)^d((2\pi)^{-2d}D_\ell)^{s/2}
\zeta_\ell(s),
\end{equation}
where $h_\ell$ is the class number and $R_\ell$ is the regulator of $\ell$, 
and
$w_\ell$ is the order of the finite group of roots of unity contained 
in $\ell.$
Using the lower bound  
$R_\ell\geqslant 0.02 w_\ell\,e^{0.1d}$ due to R.\,Zimmert [Z], we get
\begin{equation}
\frac{1}{h_{\ell,n}}\geqslant\frac1{h_\ell}\geqslant\frac{0.02}{s(s-1)}
\big( \frac{(2\pi)^s e^{0.1}}{\Gamma(s)}\big)^d
\frac1{D_\ell^{s/2}\zeta_\ell(s)}.
\end{equation}
Now from bound $(10)$ we obtain
\begin{equation}
 1>\frac{{D_\ell}^{(n-1)(n+2)/4}}
{D_k^{(n-1)/2}D_\ell^{s/2}\zeta_\ell(s)}\cdot \frac{0.02}{s(s-1)}
\big( \frac{(2\pi)^s e^{0.1}}
{\Gamma(s)}\big)^d \Big(\prod_{j=1}^{n-1}
\frac{j!}{(2\pi)^{j+1}}\Big)^{d}.
\end{equation}

\ni Letting $s = 1+\delta$, with $\delta$ in the interval $[1,10]$, and using $D_\ell\geqslant D_k^2,$ and the obvious bound $\zeta_\ell(1+\delta)\leqslant \zeta(1+\delta)^{2d},$ we get 
\begin{equation}
D_k^{1/d} \leqslant D_{\ell}^{1/2d}<\big[ \{ \frac{\Gamma(1+\delta)\zeta(1+\delta)^2}{(2\pi)^{1+\delta}
e^{0.1}}
\prod_{j=1}^{n-1}\frac{(2\pi)^{j+1}}{j!}\}\cdot
\{ 50\delta (1+\delta) \}^{1/d} \big]^{2/(n^2-2\delta-3)}.
\end{equation}

\vskip2mm

We will now prove the following simple lemma.

\begin{lemm}
For every integer $j\geqslant 2$, $\zeta_k(j)^{1/2}L_{\ell|k}(j+1)>1$.
\end{lemm}
\vskip1mm

\ni{\it Proof.} The lemma follows from the product formula for $\zeta_k(j)$ and $L_{\ell|k}(j+1)$ and the following observation.

For any positive integer $q\geqslant 2$, $$(1-\frac{1}{q^j})
(1+\frac{1}{q^{j+1}})^2 = 1-\frac{q-2}{q^{j+1}}-\frac{2q-1}{q^{2j+2}}-\frac{1}{q^{3j+2}}<1.$$
\vskip1mm

The above lemma implies that for every integer $j\geqslant 2$, $\zeta_k(j)L_{\ell|k}(j+1)>\zeta_k(j)^{1/2}>1$. Also we have the following obvious bounds for any number field $k$ of degree $d$ over $\bQ$, where, as usual, $\zeta(j)$ denotes $\zeta_{\bQ}(j)$. For every positive integer $j$, $$1< \zeta (dj)\leqslant \zeta_k(j)\leqslant \zeta(j)^d.$$ From this we obtain the following:

\vskip1mm

\begin{lemm}  Let $\cE_0=\prod^{(n-1)/2}_{j=1}\big(\zeta_k(2j)L_{\ell|k}(2j+1)\big).$  Then $\cE_0 > E_0 :=\prod_{j=1}^{(n-1)/2}\zeta(2dj)^{1/2}.$
\end{lemm}

\vskip2mm

\ni{\bf 2.11.}
To find restrictions on $n$ and $d,$ we will use the following three bounds for the relative discriminant $D_{\ell}/D_k^2$ obtained from bounds (4)-(6), 
(11), (12), and Lemma 2.

\begin{equation}
D_{\ell}/D_k^2< {\mathfrak p}_1(n,d,D_k,\delta )\ \ \ \ \ \ \ \ \ \ \ \ \end{equation} 
$$:=\big[ \frac{50\delta(1+\delta)}{E_0 D_k^{(n^2-2\delta -3)/2}}\cdot \{ 
\frac{\Gamma(1+\delta)\zeta(1+\delta)^2}{(2\pi)^{1+\delta }e^{0.1}}
\prod_{j=1}^{n-1}\frac{(2\pi)^{j+1}}{j!}\}^d \big]^{4/(n^2+ n-2\delta-4)}.$$

\begin{equation} 
D_{\ell}/D_k^2< {\mathfrak p}_2(n,d,D_k,R_{\ell}/w_{\ell},\delta) \ \ \ \end{equation}
$$:= \big[\frac{\delta(1+\delta)}{(R_\ell/w_\ell)E_0D_k^{(n^2-2\delta-3)/2}}
\cdot \{ \frac{\Gamma(1+\delta)\zeta(1+\delta)^2}{(2\pi)^{1+\delta}}\cdot
\prod_{j=1}^{n-1}\frac{(2\pi)^{j+1}}{j!}\}^d \big]^{4/(n^2+n-2\delta-4)}.$$

\begin{equation}
D_\ell/D_k^2 < {\mathfrak p}_3(n,d,D_k,h_{\ell,n})\ \ \ \ \ \ \ \ \ \ 
\end{equation}
$$:=\big[\frac{h_{\ell,n}}{E_0}\cdot \{ 
\prod_{j=1}^{n-1}\frac{(2\pi)^{j+1}}{j!}\}^d D_k^{-(n^2-1)/2}\big]^{4/(n-1)(n+2)}.$$

\ni Similarly, from bounds (4)-(6), (11), and Lemma 2 we obtain the following:
\begin{eqnarray}
&&D_k^{1/d}\leqslant D_\ell^{1/2d}<\varphi(n,d,R_{\ell}/w_{\ell},\delta)\\
&&:=\big[\{ \frac{\Gamma(1+\delta)\zeta(1+\delta)^2}{(2\pi)^{1+\delta}}
\prod_{j=1}^{n-1}\frac{(2\pi)^{j+1}}{j!}\}\cdot
\{\frac{\delta(1+\delta)}{({R_{\ell}}/{w_{\ell}}) E_0}\}^{1/d} \big]^{2/(n^2-2\delta-3)}.
\nonumber
\end{eqnarray}
\vskip3mm

\begin{center}
{\bf 3. Determining $k$}  
\end{center}
\vskip4mm

\ni{\bf 3.1.} We define $M_r(d)=\min_K D_K^{1/d},$ where the minimum is taken over all
totally real number fields $K$ of degree $d.$  Similarly, we define  $M_c(d)=\min_K D_K^{1/d},$ by taking the minimum over all
totally complex number fields $K$ of degree $d.$
 
The precise values of $M_r(d), M_c(d)$ for low values of $d$  are given in the following table (cf.\:[N]).  

$$\begin{array}{cccccccc}
d&2&3&4&5&6&7&8\\
M_r(d)^d
&5
&49
&725
&14641
&300125
&20134393
&282300416\\
M_c(d)^d&3&&117&&9747&&1257728.
\end{array}
$$

\vskip1mm
We also need the following proposition which provides lower bounds for
the discriminant of a totally real number field
in terms of its degree.

\begin{prop}
Let $k$ be a totally real number field of degree $d$, $k\neq\bQ$.  Then \\
{\rm (a)} $D_k^{1/d}\geqslant \sqrt 5>2.23.$\\
{\rm (b)} $D_k^{1/d}\geqslant 49^{1/3}>3.65$ for all $d\geqslant 3$.\\
{\rm (c)} $D_k^{1/d}\geqslant 725^{1/4}>5.18$ for all $d\geqslant 4$.\\
{\rm (d)} $D_k^{1/d}\geqslant 14641^{1/5}>6.8$ for all $d\geqslant 5$.\\
\end{prop}
\ni{\it Proof.} 
Let $g(x,d)$ and $x_0$ be as in 6.2 of [PY]. Let ${\mathfrak N}(d) = {\rm lim\,\,sup}_{x\geqslant x_0}g(x,d)$. It has been observed in [PY], Lemma 6.3, that ${\mathfrak N}(d)$ is an increasing function of $d$, and  it follows from the estimates of Odlyzko [O] that
$M_r(d)\geqslant {\mathfrak N}(d)$.  We see by a direct computation that  $g(2,9)> 9.1$. Hence, $M_r(d)\geqslant {\mathfrak N}(d)\geqslant {\mathfrak N}(9) \geqslant g(2,9) >9.1$, for all $d\geqslant 9$. 
For $1\leqslant d\leqslant 7$, from the values of
$M_r(d)$ and $M_r(d+1)$ listed above we see that  $M_r(d)\leqslant M_r(d+1).$ 
\vskip1mm

(a)--(d) now follow from the  values of $M_r(d)$, for $d\leqslant 8$, 
and
the above bound for $M_r(d)$ for $d\geqslant9.$

\bs
\ni{\bf 3.2.} We note here for latter use that except for the totally complex sextic fields with discriminants
$$-9747, \ -10051,\ -10571,\ -10816,\ -11691,\ -12167,$$
 and the totally complex quartic
fields with discriminants 
$$117,\ 125,\ 144,$$ $R_{\ell}/w_{\ell}>1/8$ 
for every number field $\ell$, see [F], Theorem B$'$.

For $r_2 = d =2,$ we have the unconditional bound 
$R_{\ell}/w_{\ell}\geqslant 0.09058$, see Theorem B$'$ and Table 3 in [F].
\vskip2mm

\ni{\bf 3.3.} For $d$ and $n$ positive integers, and $\delta \geqslant 0.02$, denote by $f(n,d,\delta)$ the expression on the extreme 
right of bounds (14) i.e.,$$ f(n,d,\delta)= \big[\{ \frac{\Gamma(1+\delta)\zeta(1+\delta)^2}{(2\pi)^{1+\delta}e^{0.1}}
\prod_{j=1}^{n-1}\frac{(2\pi)^{j+1}}{j!}\}\cdot
\{ 50\delta (1+\delta) \}^{1/d} \big]^{2/(n^2-2\delta-3)}.$$ 
For fixed $n$ and $\delta$ ($\delta\geqslant 0.02$), $f(n,d,\delta)$ clearly decreases as $d$ increases.  
\vskip1mm

We now observe that for all $n\geqslant 17$, $n!> (2\pi)^{n+1}$. From this it is easy to see that if for given $d$, $\delta$, and $n\geqslant 17$, $f(n,d,\delta) \geqslant 1$, then $f(n+1,d,\delta)< f(n, d, \delta)$, and if $f(n,d,\delta)< 1$, then $f(n+1, d, \delta)<1$. In particular, if for given $d$, and $\delta\geqslant 0.02$, $f(17, d, \delta)< c$, with $c \geqslant 1$, then $f(n,d',\delta)<c$ for all $n\geqslant 17$ and $d'\geqslant d$. 
\vskip1mm

By a direct computation we see that for $13\leqslant n \leqslant 17$, $f(n, 2,3)< 2.2$. From the  property of $f$ mentioned above, we 
conclude that $f(n,d, 3)< 2.2$ {\it for all $n\geqslant 13$, and all $d\geqslant 2$}. Now Proposition 2(a) implies that {\it for all odd} $n\geqslant 13$, $k =\bQ$.
\vskip2mm

\ni{\bf 3.4.} Now we will investigate the restriction on the degree $d$ of 
$k$ for $n\leqslant 11$ 
imposed by bound $(14).$  We get the following table by evaluating $f(n,d,\delta)$, with $n$ given in the first column, $d$ given in the second column, and $\delta$ given in the third column

$$\begin{array}{cccc}
n&d&\delta&f(n,d,\delta)<\\
11&3&2&2.6\\
9&3&1.7&3.2\\
7&4&1.5&4.1\\
5&5&1.2&6.2\\
\end{array}$$

Taking into account the upper bound in the last column of the above table, Proposition 2 implies the following:
\vskip1.5mm

\ni If $n=11,$ $d\leqslant 2.$
\vskip1mm

\ni If $n=9,$ \,\,\,$d\leqslant 2.$
\vskip1mm

\ni If $n=7,$ \,\,\,$d\leqslant 3.$
\vskip1mm

\ni If $n=5,$ \,\,\,$d\leqslant 4.$

\vskip2mm

We will now prove the following theorem by a case-by-case analysis.

\begin{theo}
If $n>7$ and the orbifold Euler-Poincar\'e characteristic of \,$\Gamma$ is
$\leqslant \chi(X_u)/n^r,$ then $d = 1$, i.e., $k = \bQ$. If $n= 7$ or $5$, and the orbifold Euler-Poincar\'e characteristic of $\Gamma$ is a submultiple of $\chi(X_u)/n^r$, then again $k =\bQ$. 
\end{theo}

\ni{\it Proof.} ({\it i}) First of all, we will show that if $n=11$, then $d$ cannot be $2$. A direct computation 
shows that $f(11,2,1.8)<2.6.$  Hence, if $n=11$ and $d=2$, then  
$D_\ell < 2.6^4< 46$. However, from the table in 3.1, we see that the smallest discriminant of 
a totally complex quartic is $117.$  So we conclude that if $n=11$, 
then $d=1.$

\vskip2mm

({\it ii}) Let us now consider the case $n=9.$
We will rule out the possibility that           
$d=2$ using bound (18).  Note that we can use the lower bound
$R_\ell/w_\ell\geqslant 0.09058$, see 3.2. We see by a direct computation that 
$\varphi(9,2,0.09058, 1.5)^4 <97.$  Hence, $D_\ell<97$ from bound (18). 
As $M_c(4)^4 =117$, $d=2$ cannot occur.  Hence, if $n=9$, then $d =1$. 
\vskip2mm

({\it iii}) We now consider the case $n=7.$
We need to rule out the possibilities that $d$ is either $3$ or $2$. 
We see from a direct 
computation that $f(7,2,1.2)<4.3$ and $f(7,3,1.4)<4.14$, where $f(n,d,\delta)$ is as in 3.3. 

   Consider first the case where $d=3$ (and $n=7$). As
$D_\ell^{1/6}<f(7,3,1.4)<4.14$, 
$D_\ell<4.14^6<5036.$  This leads to a contradiction since
according to the table in 3.1, a lower bound for the absolute value of the 
discriminant of all totally
complex sextic fields is $9747.$  Hence, it is impossible to have $d =3$ if 
$n=7$.

Consider now the case where $n=7$ and  $d=2.$  As mentioned above, 
$f(7,2,1.2)<4.3$, and hence, 
$$D_k^{1/2}\leqslant D_\ell^{1/4}< f(7,2,1.2)<4.3.$$
It follows that $D_k< 4.3^2<18.5.$  There are then the following 
five cases to discuss.
\vskip1.5mm

\ni(a) $D_k= 5$,\ \ \ $k=\bQ(\sqrt5)$\\
(b) $D_k = 8$,\ \ \ $k=\bQ(\sqrt2)$\\
(c) $D_k = 12$,\ \ $k=\bQ(\sqrt{3}$)\\
(d) $D_k =13$,\ \ $k=\bQ(\sqrt{13})$\\
(e) $D_k = 17$,\ \ $k=\bQ(\sqrt{17})$.

\vskip2mm

Case (e): We will use bound ($16$). As $R_\ell/w_\ell\geqslant0.09058$ (see 3.2), $$D_{\ell}/D_k^2<{\mathfrak p}_2(7,2,17,0.09058, 1.26)<1.1,$$ which implies that  $D_\ell = D_k^2 =17^2.$ From the table of totally complex quartics in [1], we find that there does not exist a totally complex quartic with discriminant 
$17^2$.

\vskip1mm
Case (d): $D_\ell/D_k^2< 4.3^4/13^2<2.1.$ 
Hence, $D_\ell/D_k^2=1$ or $2$. So $D_\ell=169$ or $338.$  From the table of totally complex
quartics in [1], we see that neither of these two numbers occurs as the 
discriminant
of such a field. Therefore we conclude that case (d) does not occur.
\vskip1mm

Case (c): $D_\ell/D_k^2< 4.3^4/12^2<2.4.$ 
Hence, $D_\ell/D_k^2=1$ or $2$, and $D_\ell=144$ or $288.$  Again,
from the table of totally complex
quartics in [1], we know that there is no complex quartic
with discriminant $288.$  Moreover, there is a unique totally complex 
quartic $\ell$, namely $\ell = \bQ[x]/(x^4-x^2+1)$ $=\bQ(\sqrt{-1},\sqrt{3})$, whose discriminant equals $144$. It clearly contains $k = \bQ(\sqrt{3})$. 
We will now eliminate this case using Proposition 1 (whenever we use Proposition 1 in the sequel, we will assume that the orbifold Euler-Poincar\'e characteristic of $\Gamma$ is a submultiple of $\chi(X_u)/n^r$). 

In this case, we have the following data.
$$\zeta_k(-1)=1/6, \ \ \zeta_k(-3)=23/60,\ \  
\zeta_k(-5)=1681/126,$$
$$L_{\ell|k}(-2)=1/9,\ \  L_{\ell|k}(-4)=5/3,\ \ 
L_{\ell|k}(-6)=427/3.$$

(Observe that for a positive integer $j,$ 
$\zeta_k(-(2j-1))$ and $L_{\ell|k}(-2j)$ 
are rational numbers according to well-known results of 
Siegel and Klingen. The denominators of these rational numbers can be
estimated. In this paper, we have used the software PARI together with their functional 
equations
to obtain the actual values of the Dedekind zeta and Hecke L-functions. These values have been rechecked using MAGMA. This software provides precision up to more than $40$ decimal places!)

\ni Therefore, $\mu(G(k_{v_o})/\Lambda)$ is an integral multiple of 
$$2^{-12} \zeta_k(-1)L_{\ell|k}(-2)\zeta_k(-3)L_{\ell|k}(-4)\zeta_k(-5)L_{\ell|k}(-6)=23\cdot41^2\cdot61/2^{16}\cdot3^8.$$

\ni As the numerator of this number is not a power
of $7$, according to Proposition 1 this case cannot occur.

\vskip1mm

Case (b): $D_\ell/D_k^2< 4.3^4/8^2<5.4.$
Hence, $D_\ell/D_k^2=c$ and $D_\ell=64c,$ where $c$ is a positive integer 
$\leqslant 5$. As $D_{\ell}\geqslant M_c(4)^4\geqslant 117$, the possible values of $D_{\ell}$ are $128, 192, 256, 320.$ According to the tables in [1], the only possibilities are:
\vskip1mm
\ni $D_\ell=256$: $\ell$ is obtained by adjoining a primitive $8$-th root of unity to $\bQ$; the class number of this field is $1$. 
\vskip1mm

\ni $D_\ell=320$: $\ell$ is obtained by adjoining a root of the 
polynomial $x^4-2x^3+2$ to $\bQ$, the class number of this field is also $1$.
\vskip1mm

Now,  as 
${\mathfrak p}_3(7,2,8,1)<3.1$,  from bound (17) we find that 
$D_\ell\leqslant 3\times 8^2= 192.$  So neither of the above two cases 
can occur.
\vskip1mm

Case (a): As $D_k =5$, $D_{\ell}$ is an integral multiple of $25$. We will now use bound ($16$) to find an upper bound for $D_\ell/D_k^2$,
making use of the estimate of Friedman [F] mentioned in 3.2
that $R_\ell/w_\ell >1/8$ if $D_\ell\neq 125.$ We find that 
$D_\ell/D_k^2<{\mathfrak p}_2(7,2,5,1/8,1.3)<8.7.$ 
So $D_\ell =25c$, where c is a positive integer 
$\leqslant 8$.  Since the smallest discriminant of
a totally complex quartic
is $117$, $c\geqslant 5.$  
Hence, $5\leqslant c\leqslant 8.$  
The possible values of $D_{\ell}$ are therefore $125,\ 150,\ 175,\ 200.$  From the tables in [1] we see that there is no totally complex quartic field with discriminant $150$, $175$ or $200$, whereas the field $\ell$ 
obtained by adjoining a primitive 
$5$th root of unity to $\bQ$ is the unique totally complex quartic field with $D_{\ell} = 125$. It is a cyclic extension of $\bQ$, and it contains 
$k = \bQ(\sqrt5)$.  We will use Proposition 1 to eliminate this case. 
In this case, we have the following data.
$$\zeta_k(-1)=1/30, \ \ \zeta_k(-3)=1/60,\ \  \zeta_k(-5)=67/630,$$
$$L_{\ell|k}(-2)=4/5,\ \  L_{\ell|k}(-4)=1172/25,\ \ 
L_{\ell|k}(-6)=84676/5.$$
Hence $\mu(G(k_{v_o})/\Lambda)$ is an integral multiple of 
$$2^{-12} \zeta_k(-1)L_{\ell|k}(-2)\zeta_k(-3)L_{\ell|k}(-4)\zeta_k(-5)L_{\ell|k}(-6)= {67\cdot 293\cdot 21169}/{2^{10}\cdot 3^4 \cdot 5^7\cdot 7}.$$
Again, as the numerator of this number is not a power
of $7$, according to Proposition 1 this case cannot occur. 
\vskip2mm

({\it iv}) Finally we take-up the case $n=5.$  We will rule out the possibilities that $d=4$, $3$ or $2$.
\vskip2mm

({\bf 1}) Consider first the case where $n =5$ and $d=4.$  Bound  (14) with $\delta=1$
leads to $D_{\ell}^{1/8}<f(5,4,1)<6.4.$ 
Now from Table 2 of [F] we find that  
$R_\ell/w_\ell\geqslant 0.1482$. Next we use bound (18) to conclude that $D_k^{1/4}\leqslant D_{\ell}^{1/8}<{\varphi}(5,4,0.1482,1.2)<6.05.$
As $6.05^4< 1340,$ $D_k < 1340.$  From the list of quartics with small discriminants given in [1], we see that the only integers smaller than $1340$ which are the discriminant of a totally
real quartic $k$ are $725$ and $1125.$ 
Moreover, for either of these two integers, there is a unique totally real quartic field $k$ whose discriminant 
is that integer. Each of these fields has class number $1$. 
\vskip1mm

If $D_k=1125,$ 
$$D_\ell/D_k^2< {6.05^8}/{1125^2}<2.$$
So $D_\ell/D_k^2=1.$
 This implies
that 
$D_\ell=1125^2=1265625.$
\vskip1mm

If $D_k=725,$  
$$D_\ell/D_k^2< {6.05^8}/{725^2}<4.$$
Hence $D_\ell=725^2 c$ with  $c\leqslant3$.  In particular,
$D_\ell\leqslant 1576875.$
\vskip1mm
At our request, Gunter Malle has shown by explicit computation{\footnote {Malle used the following procedure in his computation. Any quadratic extension of $k$ is of the form $k(\sqrt{\alpha})$, with $\alpha$ in the ring of integers ${\mathfrak o}_k$ of $k$. As the class number of any totally real quartic $k$ presently under consideration is $1$, ${\mathfrak o}_k$ is a unique factorization 
domain. Now using factorization of small primes and explicit generators of the group of units of $k$, he listed all possible $\alpha$ modulo squares, and then for each of the $\alpha$, the discriminant of $k(\sqrt{\alpha})$ could be 
computed.}} that there is exactly one pair 
of number fields $(k,\ell)$ with $(D_k, D_{\ell})$ 
among the four possiblities above.  $k$ (resp.,\,\,$\ell$) is obtained by adjoining a root of $x^4-x^3-4x^2+4x+1$ (resp.,\,\,a primitive $15$th root of unity which is a root of $x^8-x^7+x^5-x^4+x^3-x+1$) to $\bQ$. For this pair $D_k = 1125$, $D_{\ell} = 1125^2 = 1265625$, and the class number of $\ell$ is $1$. We will now employ Proposition 1 to eliminate this case.   We have the following values of the Dedekind zeta and Dirichlet $L$-functions for this pair 
$(k,\ell)$.
$$ \zeta_k(-1)=4/15,\ \zeta_k(-3)=2522/15,\ L_{\ell|k}(-2)=128/45, \ 
L_{\ell|k}(-4)=2325248/75.$$ From which we conclude that $\mu(G(k_{v_o})/\Lambda)$ is an integral 
multiple of 
$$2^{-16} \zeta_k(-1)L_{\ell|k}(-2)\zeta_k(-3)L_{\ell|k}(-4)
={2^2\cdot13\cdot31\cdot97\cdot293}/{3^5\cdot5^5}.$$
As the numerator of this number is not a power of $5$, Proposition 1 rules out this case.

\vskip2mm

({\bf 2}) We will consider now the case where $n=5$ and $d=3.$  As $\ell$ is a totally complex 
sextic 
field, from 3.2 we know that  $R_\ell/w_\ell >1/8$ unless $\ell$ is a totally complex sextic field whose discriminant equals one of the six negative 
integers listed in 3.2.  Now using this lower bound for $R_{\ell}/w_{\ell}$, we deduce from (18) that 
$D_k \leqslant D_\ell^{1/2}<{\varphi}(5,3,1/8,1)^3<6.24^3 <243$.   On the other hand, if $\ell$ is a totally complex sextic field whose discriminant equals one of the six negative integers listed in 3.2, then $D_k\leqslant 12167^{1/2}<111$. Hence, if $n=5$, $d=3$, then 
$D_k<243$. From Table B.4 in [Co] of discriminants  
of totally real cubic number fields we infer that $D_k$ must equal one of the 
following five integers: $49,\ 81,\ 148,\ 169,$ and $229.$
\vskip1mm

\ni$\bullet$ If $D_k=229,$ $D_\ell/D_k^2< 6.24^6/229^2<1.2.$
Hence,  $D_\ell=229^2=52441.$
There are however no such totally complex sextic fields according to [1].
\vskip1mm

\ni$\bullet$ If $D_k = 169$ or $148$, $D_{\ell}\geqslant
D_k^2\geqslant 148^2> 12167$, and hence $R_\ell/w_\ell >1/8$, see
3.2. We will now use bound (16). As ${\mathfrak p}_2(5,3,169,1/8,1.1)<1.9$, and ${\mathfrak p}_2(5,3, 148, 1/8, 1.1)<2.3$, $D_{\ell}$ must equal $cD_k^2$ for some $c\leqslant 2$. 
None of these 
appear in 
 the table t60.001 of totally complex sextics in [1].

\vskip1mm

\ni$\bullet$ If $D_k=81,$ then $81^2|D_{\ell}$, but none of the six negative integers listed in 3.2 are divisible by $81^2$. Hence, $R_{\ell}/w_{\ell}> 1/8$. 
We will again use bound (16). We see by a direct computation that 
${\mathfrak p}_2(5,3,81,1/8,1.1)<6.2$. Therefore, if $D_k = 81$, then $D_{\ell}= cD^2_k$, with  $1\leqslant c\leqslant 6$. But from the table t60.001 in [1] we see that there is no totally complex 
sextic field the absolute value of whose discriminant equals $81^2 c$, with $1\leqslant c\leqslant 6$, except for $c =3$. Thus   
$D_\ell$ can only be $3\times81^2= 19683$.
\vskip1mm

Let $k$ be the field obtained by adjoining a root of $x^3-3x-1$ to $\bQ$, and 
$\ell$ its totally complex quadratic extension obtained by adjoining a primitive $9$th root of unity to $\bQ$. Then $k$ (resp.,\,\,$\ell$) is the unique totally real 
cubic (resp.,\,\,totally complex sextic) field with $D_k=81$ (resp.,\,\,$D_{\ell} = 19683$). In this case, we have the following data on the values of the zeta and $L$-functions.
$$ \zeta_k(-1)=-1/9,\ \zeta(-3)=199/90,\  L_{\ell|k}(-2)=-104/27,
\ L_{\ell|k}(-4)=57608/9.$$ From which we conclude that $\mu(G(k_{v_o})/\Lambda)$ is an integral multiple of 

$$2^{-12} \zeta_k(-1)L_{\ell|k}(-2)\zeta_k(-3)L_{\ell|k}(-4)={13\cdot 19\cdot 199\cdot379}/{2^7\cdot3^9\cdot 5}.$$
As the numerator of this rational number is not a power of $5$, according to Proposition 1 this case cannot occur.

\vskip2mm

\ni$\bullet$ If $D_k=49,$ then $D_{\ell}$ is divisible by $49^2$, but none of the six negative integers listed in 3.2 are divisible by $49^2$. So $R_{\ell}/w_{\ell}>1/8$. We apply bound (16) to obtain
$D_\ell/D_k^2< {\mathfrak p}_2(5,3,49,1/8,1.2)$ $<14.3.$ Hence, $D_\ell=49^2 c$, with $1\leqslant c\leqslant 14$. On the other hand, the table in 3.1 
implies that
$c>9747/49^2>4.$  Therefore, we need only consider $5\leqslant c\leqslant  14$. From the table t60.001 in [1] we see that among these ten integers, 
$7\times 49^2= 
16807$ is the only one which equals $D_{\ell}$ of a totally complex sextic $\ell$. This $\ell$ is obtained by adjoining a primitive $7$th root of unity to $\bQ$ and it contains the totally real cubic field $k$ obtained by adjoining a root of $x^3-x^2-2x+1$ to $\bQ$. 
It is easy to see that $D_k=49$ in this case. We have the following data on the values of the zeta and $L$-functions.
$$ \zeta_k(-1)=-1/21,\ \zeta(-3)=79/210,\  L_{\ell|k}(-2)=-64/7,
\ L_{\ell|k}(-4)=211328/7.$$ From which we conclude that $\mu(G(k_{v_o})/\Lambda)$ is an integral multiple of 

$$2^{-12} \zeta_k(-1)L_{\ell|k}(-2)\zeta_k(-3)L_{\ell|k}(-4)={13\cdot 79\cdot 127}/{3^2\cdot 5\cdot7^4}.$$

\ni Again, as the numerator of this rational number is not a power of $5$, according to Proposition 1 this case cannot occur.

\vskip2mm

({\bf 3}) We will consider now the case $n=5$, $d=2.$ We recall the lower 
bound $R_\ell/w_\ell \geqslant 0.09058$ from 3.2. From bound (18) we obtain that 
$D_k^{1/2}\leqslant D_\ell^{1/4}< {\varphi}(5,2,0.09058,1)< 6.7$. Since 
$6.7^2<45$, $D_k\leqslant 44.$ 
It follows that the discriminant $D_k$ of the real quadratic field $k$ can
only be one of the following integers,
$$5,\ 8,\ 12,\ 13,\ 17,\ 21,\ 24,\ 28,\ 29,\ 33,\ 37,\ 40,\ 41,\ 44.$$

\vskip2mm

\ni$\bullet$ If $D_k\geqslant 37$, then $D_\ell/D_k^2< 6.7^4/37^2<2.$  
In these cases, $D_\ell=D_k^2$ is one of the following integers 
$1369,\ 1600,\ 1681,\ 1936.$  Of these only $1600$ and $1936$ appear as 
the discriminant of
a totally complex quartic $\ell$, check [1]. Moreover, there is a unique 
totally complex quartic $\ell$ with $D_{\ell}= 1600$ (resp.,\,\,$D_{\ell}= 
1936$). The class number of both of these quartics is $1$. Now we will use 
bound (17). Since ${\mathfrak p}_3(5,2,40,1)<0.6<1$ and  
${\mathfrak p}_3(5,2, 44,1)<0.5<1,$ if either $D_k =40$ or $44$, then $D_{\ell}/D_k^2<1$, which is impossible.

\vskip1.5mm

\ni$\bullet$ If $D_k=33,$ then $D_{\ell}\geqslant 33^2= 1089$, and hence 
$R_\ell/w_\ell>1/8$, see 3.2. Now  from bound (16) we obtain that 
$D_\ell/D_k^2< {\mathfrak p}_2(5,2,33,1/8,1)<2.$ Hence, 
$D_\ell=D_k^2=1089$. There is a unique totally
complex quartic $\ell$ whose discriminant is $1089$. Its class number is $1$. 
Now we apply bound (17), $1\leqslant D_\ell/D_k^2< {\mathfrak p}_3(5,2,33,1)<0.77$, to reach a contradiction.  

\vskip1.5mm

\ni$\bullet$ If $D_k=29,$ then $D_\ell/D_k^2< 6.7^4/29^2<3.$
Hence, $D_\ell/D_k^2=1$ or $2$ . Therefore, $D_\ell=29^2=841$ or $1682$. Neither  
of these two integers is the discriminant of a totally complex quartic ([1]). 

\vskip1.5mm

\ni$\bullet$ If $D_k=17$ or $13$, then $D_{\ell}\geqslant 169$, and hence $R_{\ell}/w_{\ell}>1/8$ from 3.2. Now we will use bound 
(16). As ${\mathfrak p}_2(5,2,17,1/8,1)<4.7$, and 
${\mathfrak p}_2(5,2,13,1/8,1)<7.2$, 
$D_\ell =17^2 c$, with $1\leqslant c\leqslant 4$, or $D_{\ell} = 13^2 c$, 
with $1\leqslant c\leqslant 7$. But of these eleven integers none
 appears as the discriminant of a totally complex 
quartic field. 
\vskip1.5mm

\ni$\bullet$ To eliminate the remaining cases (namely, where $D_k = 5,\ 8,\ 12,\ 21,\ 24$ or $28$), we will use Proposition 1. Let us assume in the rest of this section that $D_k$ is one of the following six integers: $5$, $8$, $12$, $21$, $24$, $28$. As $D_{\ell}$ is an integral multiple of $D_k^2$, we conclude from 3.2 that 
unless $D_{\ell} = 125$ or $144$, $R_{\ell}/w_{\ell}>1/8$. We will now use upper bounds (16) and (17) for $D_{\ell}/D_k^2$ to make a list of the pairs $(k,\ell)$ which can occur.
\vskip1.5mm

\ni (i) As ${\mathfrak p}_2(5,2,28,1/8,1)<2.1$, if $D_k = 28$, then
$D_{\ell}= 28^2 c$, with $c=1$ or $2$. We see from [1] that the class number of any totally complex quartic $\ell$ with $D_{\ell}= 28^2$ or $2\times 28^2$ is $1$. Now we note that ${\mathfrak p}_3(5,2,28,1)<1.1$. Hence $D_{\ell}$ can only be $28^2 = 784$. The corresponding quartic field is $\ell = \bQ[x]/(x^4-3x^2+4)=\bQ(\sqrt{-1},\sqrt{7})$, which contains $k = \bQ(\sqrt{7})$. We shall denote this pair $(k, \ell)$ by ${\mathfrak C}_1$.  
\vskip1.5mm

\ni (ii) As ${\mathfrak p}_2(5,2,24,1/8,1)<2.6$, if $D_k = 24$, then $D_{\ell}
= 24^2 c$, with $1\leqslant c\leqslant 2$. Of these integers, only $24^2 =576$ is the discriminant of a totally complex quartic. 
There are two totally complex quartics with discriminant $576$, but only one of them contains $k = \bQ(\sqrt{6})$. This quartic is $\ell = \bQ[x]/(x^4-2x^2+4)=\bQ(\sqrt{-3},\sqrt{6})$. We shall denote this pair $(k,\ell )$ by ${\mathfrak C}_2$.    
\vskip1.5mm

\ni (iii) As ${\mathfrak p}_2(5,2,21,1/8,1)<3.3$, if $D_k = 21$, 
then $D_{\ell}= 21^2 c$, with $1\leqslant c\leqslant 3$. Of these three integers, 
only $21^2=441$ is the discriminant of a totally complex quartic $\ell$. This quartic is $\ell =\bQ[x]/(x^4-x^3-x^2-2x+4)=\bQ(\sqrt{-3},\sqrt{-7})$, and it clearly contains $k =\bQ(\sqrt{21})$. We shall denote this pair $(k,\ell)$ by ${\mathfrak C}_3$. 
\vskip1.5mm

\ni (iv) As ${\mathfrak p}_2(5,2,12,1/8,1)<8.3$, if $D_k =12$, then $D_{\ell} = 12^2 c$, with $1\leqslant c \leqslant 8$. Among these, only for $c=1,\ 3,\ 4,$ and $7$, there exists a totally complex quartic $\ell$ with $D_{\ell} = 12^2 c$, and all these quartics have the class number $1$. Now we note that ${\mathfrak p}_3(5,2,12, 1)<4.4$, which implies that $c\leqslant 4$; 
i.e., $c= 1,\,3$, or $4$. The quartics corresponding to $c = 3$ and $4$ do not contain $\bQ(\sqrt{3})$. 
As we observed while dealing with Case (c) in ({\it iii}) above, there is a unique totally complex quartic $\ell$, namely $\ell = \bQ[x]/(x^4-x^2+1)$ $=\bQ(\sqrt{-1},\sqrt{3})$, whose discriminant equals $12^2 =144$. It contains $k = \bQ(\sqrt{3})$. The pair $\big(\bQ(\sqrt{3}),
\bQ(\sqrt{-1},\sqrt{3})\big)$ will be denoted by ${\mathfrak C}_4$.    
\vskip1.5mm

\ni (v) As ${\mathfrak p}_2(5,2,8,1/8,1)< 16.2$, if $D_k = 8$, then 
$D_{\ell}= 8^2 c$, with $1\leqslant c \leqslant 16$. Among these, only for 
$c=4, \  5,\ 8, \ 9,$ and $13$, there exists a totally complex 
quartic field with discriminant $8^2 c$, and all these quartics have the class number $1$. Now
we observe that ${\mathfrak p}_3 (5,2,8,1)<8.7$, which implies that
$c= 4,\ 5$ or $8$. 
There is only one
 totally complex quartic field $\ell$ containing $k = \bQ(\sqrt{2})$, 
with discriminant as above. This is 
$\ell = \bQ[x]/(x^4+1)=\bQ(\sqrt{-1},\sqrt{2})$  (with $D_{\ell} = 256$). The corresponding pair 
$(k,\ell)$ will be denoted by ${\mathfrak C}_5$.       
\vskip1.5mm

\ni (vi) As ${\mathfrak p}_2(5,2,5,1/8, 1)<35.5$, and $D_{\ell}\geqslant 117$, see 3.1, if $D_k =5$, 
then $D_{\ell} = 25c$, with $5\leqslant c \leqslant 35$. Among these, only for $c= 5,\ 9$ and $16$, there exists a totally complex quartic field with discriminant $25c$.  Thus the possible values of $D_{\ell}$ are $125$, $225$ and $400$. There are precisely three totally complex quartic fields containing $k =\bQ(\sqrt{5})$ and with discriminant in $\{125,\,225,\,400\}$. These are $\ell =\bQ[x]/(x^4-x^3+x^2-x+1)$ (= the field obtained by adjoining a primitive $5$th root of unity to $\bQ$, its discriminant is $125$), $\ell = \bQ[x]/(x^4-x^3+2x^2+x+1)=\bQ(\sqrt{-3},\sqrt{5})$ (with discriminant $225$), and $\ell = \bQ[x]/(x^4+3x^2+1)=\bQ(\sqrt{-1},\sqrt{5})$ (with discriminant $400$). The corresponding pairs $(k,\ell)$ will be 
 denoted by ${\mathfrak C}_6$, ${\mathfrak C}_7$ and ${\mathfrak C}_8$ respectively.    
\vskip2mm

We observe that in all the above cases, the conclusion of 
Proposition 1 is violated, see the last column of the table below, where $\cR= 2^{-8}\zeta_k(-1)L_{\ell |k}(-2)\zeta_k(-3)L_{\ell |k}(-4)$ is as in (7) for $n=5$ and $d=2$.  Hence
none of these cases can occur. We have thus completely proved Theorem 1. 
\vskip2.5mm

$$\begin{array}{ccccccccc}
(k,\ell)&\ \ \zeta_k(-1)&\ \ \zeta_k(-3)&\ \ L_{\ell|k}(-2)&\ \ L_{\ell|k}(-4)&{\cR}\\
{\mathfrak C}_1&2/3&113/15&8/7&80&113/3^2\cdot 7\\
{\mathfrak C}_2&1/2&87/20&2/3&38&19\cdot29/2^9\cdot5\\
{\mathfrak C}_3&1/3&77/30&32/63&64/3&2^2\cdot 11/3^{5}\cdot 5\\
{\mathfrak C}_4&1/6&23/60&1/9&5/3&23/2^{11}\cdot3^5\\
{\mathfrak C}_5&1/12&11/120&3/2&285/2&11\cdot19/2^{15}\\
{\mathfrak C}_6&1/30&1/60&4/5&1172/25&293/2^7\cdot3^2\cdot5^5\\
{\mathfrak C}_7&1/30&1/60&32/9&1984/3&31/3^5\cdot5^2\\
{\mathfrak C}_8&1/30&1/60&15&8805
&587/2^{11}.
\end{array}
$$

\bs
\ni
\begin{center}
{\bf 4. Restrictions on $\ell$ and the main result}  
\end{center}
\vskip4mm

\ni{\bf 4.1.} We shall assume in the sequel that $k = \bQ$. (We have proved in the preceding section that this is the case if $n>7$, or if $n=7$ or $5$ and the orbifold Euler-Poincar\'e characteristic of $\Gamma$ is a submultiple of $\chi(X_u)/n^r$.) 
Then  $\ell=\bQ(\sqrt{-a})$ for some square-free positive integer $a.$
By setting $d = 1$ and $D_k =1$ in bound (13) we obtain
$$1 > \frac{D_\ell^{(n-1)(n+2)/4}}{D_{\ell}^{s/2}\zeta_{\ell}(s)}
\cdot\frac{0.02}{s(s-1)}\cdot\frac{(2\pi)^s e^{0.1}}{\Gamma(s)}\cdot
\prod_{j=1}^{n-1}\frac{j!}{(2\pi)^{j+1}}.$$

Using the obvious bound $\zeta_{\ell}(s)\leqslant \zeta(s)^2$, and by 
setting $s=1+\delta,$ we derive from the above that  
\begin{equation}
 D_\ell < \{50\delta(1+\delta)e^{-0.1}\Gamma(1+\delta)(2\pi)^{-1-\delta}\zeta(1+\delta)^2\prod_{j=1}^{n-1}\frac{(2\pi)^{j+1}}{j!}\}^{4/(n^2+n-2\delta-4)}.
\end{equation}

\ni{\bf 4.2.}
Denote by $\mathfrak{d} (n,\delta)$ the right hand side of the above bound.  We see, as in 3.3, that for a fixed value of $\delta$, $\mathfrak{d} (n,\delta)$ decreases as $n$ increases provided $n\geqslant 19$. We obtain the following
upper bound for 
$\mathfrak{d} (n,\delta)$ for $n$ listed in the first column and $\delta$ listed in the second column of the following table:
$$\begin{array}{ccc}
n&\delta&D_{\ell}<\mathfrak{d} (n,\delta)<\\
19&2&2.2\\
17&2&2.7\\
15&2&3.4\\
13&2&4.5\\
11&2&6.2\\
9&2&9.4\\
7&1&15.7\\
5&0.5&37.4
\end{array}$$

The bound for $D_{\ell}$ given by the bound for $\mathfrak{d}(n,\delta)$ in the 
above table restricts the possibilities for $n$ and $\ell.$  In particular, 
since an imaginary quadratic field has discriminant at least
$3$, we deduce
from the above table and the monotonicity of $\mathfrak{d} (n,\delta)$ for
a fixed $\delta$ that {\it it is impossible for $n$ to be larger than $15$}. We recall that for $\ell = \bQ (\sqrt{-a} )$, where $a$ is a 
square-free positive integer, $D_{\ell} =a$ if $a\equiv 3$ (${\rm mod}\ 4$), and $D_{\ell} =4a$ otherwise. From the above table we now obtain the following enumeration of  all possible $n$ and $\ell$.
\vskip2mm

\ni{\rm (a)} $n\leqslant 15.$ 
\vskip1mm

\ni{\rm (b)} The number field $\ell$ equals $\bQ(\sqrt{-a})$, where for all odd $n$, $5\leqslant n\leqslant 15$,  the possible values of $a$ are listed below:

$$\begin{array}{cc}
n&a\\
15&3\\
13&1,3\\
11&1,3\\
9&1,2,3,7\\
7&1,2,3,7,11,15\\
5&1,2,3,5,6,7,11,15,19,23,31,35.
\end{array}$$

\vskip2mm
\ni{\bf 4.3.}
It is known that the class number of the fields $\ell$ appearing in the above 
table is $1$, except when $a=5$, $6$, $15$, or $35$, in which cases $\ell$ 
has the class number $2$, or $a=23$, or $31$, in which cases $\ell$ has the  
class number $3$.  Hence from (10) we get the following bound:
\begin{equation}
D_\ell < [h_{\ell,n}\prod_{j=1}^{n-1}\frac{(2\pi)^{j+1}}{j!}]^{4/(n-1)(n+2)},
\end{equation}
where $h_{\ell,n}$ can only be $1$ or $3$ since $n$ is odd. Let $\lambda (n,h_{\ell,n})$ be the function on the right hand side of the above bound.
Direct computation yields the following table.
$$\begin{array}{ccccccc}
n&15&13&11&9&7&5\\
\lambda(n,3)<&3.3&{\ \ }&{\ \ }&8.1&{\ \ }&{\ \ }\\
\lambda(n,1)<&3.3&4.2&5.5&7.7&11.2&17.6.
\end{array}$$

\vskip2mm

Using the above table, and upper bound (20) for $D_{\ell}$, we 
conclude the following.

\begin{prop}
The only possibilities for the number field $\ell=\bQ(\sqrt{-a})$ are those listed in the following table.\\
$$\begin{array}{cc}
n&a\\
15&3\\
13&1,3\\
11&1,3\\
9&1,3,7\\
7&1,2,3,7,11\\
5&1,2,3,7,11,15.
\end{array}$$

\end{prop}

\vskip1mm

\ni {\bf 4.4.} In the considerations so far we did not need to assume that $\Gamma$ is cocompact. {\it We will henceforth assume that $\Gamma$ is cocompact}, and make use of the description of $G$ given in the introduction. Let $\ell$, the division algebra 
$\cD$, and the hermitian form $h$ be as in there.

If $\cD = \ell$, then $h$ is an hermitian form on $\ell^n$ and its signature over $\bR$ is $(n-m,m)$, $n>m>0$. The hermitian form $h$ gives us a quadratic form $q$ on the $2n$-dimensional $\bQ$-vector space $V= \ell^n$ defined as follows:  $$q(v)= h(v,v) \ \ \ {\rm{for}}\ \ v\in V.$$
The quadratic form $q$ is isotropic over $\bR$, and hence by Meyer's theorem it is isotropic over $\bQ$ (cf.\,[Se2]). This implies that $h$ is isotropic, and hence so is $G/\bQ$. Then by Godement's compactness criterion, $\Gamma$ is 
noncocompact, which is contrary to our hypothesis. We conclude therefore 
that $\cD\ne \ell$, and so it is a nontrivial central simple division 
algebra over $\ell$. 
\vskip1mm

    From the classification of central simple division algebras over $\ell$, which admit an involution of the second kind, we know that the set $\cT_0$ of rational primes $p$ which split in $\ell$, but the group $G$ does not split over 
$\bQ_p$, is nonempty. Note that $\cT_0\subset \cT$, where $\cT$ is as 
in 2.4, and $p\in \cT_0$ if, and only if, ${\bQ}_p\otimes_{\bQ}\cD=({\bQ}_p\otimes_{\bQ}\ell )\otimes_{\ell}\cD$ is isomorphic to $M_r({\mathfrak D}_p)\oplus 
M_r({\mathfrak D}_{p}^{o})$, where ${\mathfrak D}_p$ is a noncommutative central division algebra over 
$\bQ_p$, ${\mathfrak D}_p^o$ is its opposite, and $r$ is a positive integer. We shall denote the degree of ${\mathfrak D}_p$ by $d_p$ in the sequel.

\vskip1mm
    
\ni{\bf 4.5.}  Now we will use the Euler product $\cE$ appearing in the volume formula (5) to eliminate all but the pair  $(n,a) = (5,7)$ appearing in 
Proposition 3. Recall from 2.7 that  
\begin{eqnarray*}
\cE&=&\prod_{p\in\cT\cup\cT'}e^{\prime}(P_p)
\prod_{j=1}^{(n-1)/2}\big(\zeta(2j)L_{\ell|\bQ}(2j+1)\big)\\
&=&\cE_1\cE_2\cE_3,
\end{eqnarray*}
where
\begin{eqnarray*}
\cE_1&=&\prod_{p\in\cT\cup\cT'}e^{\prime}(P_p),\\
\cE_2&=&\prod_{j=1}^{(n-1)/2}\zeta_{\ell}(2j+1),\\
\cE_3&=&\prod_{j=1}^{(n-1)/2}\frac{\zeta(2j)}{\zeta(2j+1)}.
\end{eqnarray*}
In the above we have used the simple fact that 
$L_{\ell|\bQ}(j)=\zeta_{\ell}(j)/\zeta(j).$
\vskip1mm

\ni{\bf 4.6.}
Clearly, $\cE_2>1$ since each factor in the product formula for 
$\zeta_\ell(2j+1)$, for $j>0$, is greater than $1$. Also, $e'(P_p)$ is an integer for all $p$, and for $p\in\cT$, $e'(P_p)>n$ (see 2.1 and 2.3). Now from (4), (5), (6) and
(9) we obtain
\begin{eqnarray*}
D_\ell&\leqslant& \big(h_{\ell,n}\frac{n^{\#\cT}}{\cE}\prod_{j=1}^{n-1}\frac{(2\pi)^{j+1}}{j!}\big)^{4/(n-1)(n+2)}\\
&<&\big(h_{\ell,n}\frac{n^{\#\cT}}{\cE_1\cE_3}\prod_{j=1}^{n-1}\frac{(2\pi)^{j+1}}{j!}\big)^{4/(n-1)(n+2)}\\
&\leqslant&\Big(h_{\ell,n}\cdot\prod_{p\in\cT_0}\frac{n}{e'(P_p)}\cdot\prod_{j=1}^{(n-1)/2}\frac{\zeta(2j+1)}{\zeta(2j)}
\cdot\prod_{j=1}^{n-1}\frac{(2\pi)^{j+1}}{j!}\Big)^{4/(n-1)(n+2)}.
\end{eqnarray*}

It follows from 2.1 and 2.3(ii) that for $p\in \cT_0$, $e'(P_p)$ is an integral
multiple of $$\frac{\prod_{j=1}^{n}(p^j -1)}{\prod_{j=1}^{n/d_p}(p^{jd_p}-1)},$$
where $d_p>1$ and $d_p|n.$  Let $q$ be the largest prime belonging 
to $\cT_0$. Then 
$$\prod_{p\in\cT_0}\frac{e'(P_p)}{n}\geqslant\frac{1}{n}\cdot \frac{\prod_{j=1}^{n}(q^j -1)}{\prod_{j=1}^{n/d_{q}}(q^{jd_{q}}-1)},$$
which implies that
$$D_\ell < L(n,d_q,q,h_{\ell,n}),$$
where for any divisor $d$ of $n$, $$L(n,d,q,h_{\ell,n})= 
\Big(nh_{\ell,n}\cdot\frac{\prod_{j=1}^{n/d}(q^{jd}-1)}
{\prod_{j=1}^{n}(q^j -1)}
\cdot\prod_{j=1}^{(n-1)/2}\frac{\zeta(2j+1)}{\zeta(2j)}
\cdot\prod_{j=1}^{n-1}\frac{(2\pi)^{j+1}}{j!}\Big)^{4/(n-1)(n+2)}.$$

 Note that $L(n,d, q,h_{\ell,n})$ is decreasing in $q$ if the other 
three arguments are fixed. Also note that $L(n, d_{q},q,h_{\ell,n})\leqslant L(n,d,q,h_{\ell,n})$, where $d$ is any divisor of $d_{q}$. 
\vskip1mm

 Let $a$ be a square-free positive integer. We recall now the following 
well-known fact (cf.\,[BS]).
\begin{lemm}
(a) An odd prime $p$ splits in $\ell 
= {\bQ (\sqrt{-a})}$ if, and only if, $p$ does not divide $a$, and $-a$ is a 
square modulo $p.$\\
(b) $2$ splits in $\ell$ if, and only if,  
$a \equiv -1 \ ({\rm mod}\ 8)$.\\
(c) A prime $p$ ramifies in $\ell$ if, and only if, $p|D_{\ell}.$
\end{lemm}

As $q\in \cT_0$, $q$ splits in 
$\ell.$ Thus if $p =p_a$ is the smallest prime splitting in 
$\ell = \bQ(\sqrt{-a})$, then $q\geqslant p$. Hence, $D_{\ell}< L(n,d_{q}, q, h_{\ell,n})\leqslant L(n,d_{q}, p, h_{\ell,n})$. 
\vskip1mm

We easily see using Lemma 3 that the smallest prime splitting in 
$\ell = \bQ(\sqrt{-a})$ for $a = 1$, $2$, $3$, $7$, $11$ and $15$ are respectively $5$, $3$, $7$, $2$, $3$ and $2$. The class number $h_{\ell}$ of $\ell=\bQ(\sqrt{-a})$, for $a = 1$, $2$, $3$, $7$, $11$ and $15$ are $1$, $1$, $1$, $1$, $1$ and $2$ respectively. Now we see by a simple 
computation that for the pairs $(n,a)$ appearing in Proposition 4, 
$L(n,d, p, h_{\ell,n})< D_{\ell}$, for any prime divisor $d$ of $n$, except for 
$(n,a) = (5,7)$. Moreover, $L(5,5,2,1)> D_{\bQ(\sqrt{-7})} = 7$, but for any $q>2$, $L(5,5,q,1)< 7$. We conclude therefore the following.
\begin{theo}
The only possibilities for $\ell$, $n$ and $\cT_0$ are $\ell = \bQ(\sqrt{-7})$, $n = 5$ and $\cT_0 = \{ 2\}$. 

In particular, ${\rm PU}(n-m,m)$, with $n$ odd, and  $0<m<n$,  can contain a cocompact arithmetic subgroup 
whose orbifold Euler-Poincar\'e characteristic is $\chi(X_u)/n$, where $X_u$ is the compact dual of the symmetric space of ${\rm PU}(n-m,m)$,  only if $n =3$ or $5$.  
\end{theo}         

\vskip4mm

\begin{center}
\ni {\bf 5. Four arithmetic fake $\bP_{\bC}^4$ and four arithmetic fake ${\bf Gr}_{2,5}$}
\end{center}
\vskip4mm

\ni{\bf 5.1.} Let now $k=\bQ$, $\ell = \bQ(\sqrt{-7})$, and $\cD$ be a division algebra with center $\ell$ and of degree $5$ such that for every place $v$ of $\ell$ 
not lying over $2$, $\ell_v\otimes_{\ell}\cD$ is the matrix algebra $M_5(\ell_v )$, and the invariant of $\cD$ at $v'$ is $a/5$ and at $v''$ it is $-a/5$, where $v'$ and $v''$ are the places of $\ell$ lying over $2$, and $a$ is a positive integer 
less than $5$. Let $m= 1$ or $2$. Then $\cD$ admits an involution $\sigma$ of the second kind such that if $G$ is the simply connected simple algebraic $\bQ$-group with $$G(\bQ) = \{x\in \cD^{\times}\ |\ x\sigma(x)=1\ \ {\rm{and}}\ \ {\rm{Nrd}}\, x = 1\},$$ then $G(\bR)$ is isomorphic to ${\rm SU}(5-m,m)$. We note that by varying $a$, and for a given $m$,  varying the involution $\sigma$ of $\cD$, we get exactly {\it two} distinct 
groups $G$ up to $\bQ$-isomorphism. 
\vskip1mm

Let $G$ be as above. We fix a maximal compact-open subgroup $P= \prod P_q$ of the group $G(A_f)$ of finite ad\`eles of $G$, where for all $q\ne 2$, $7$, $P_q$ is a hyperspecial parahoric subgroup of $G(\bQ_q)$, $P_2 = G(\bQ_2)$, and $P_7$ is 
a {\it special} maximal parahoric subgroup of $G(\bQ_7)$ (we note that there 
are exactly {\it two} such parahoric subgroups containing a given Iwahori subgroup 
of $G(\bQ_7)$ and they are nonisomorphic as topological groups, cf.\:[T2]). Let $\Lambda = G(\bQ)\cap P$. Then $\Lambda$, considered as a subgroup of $G(\bR)$, is a principal arithmetic subgroup. The following lemma implies that 
$\Lambda$ {\it is torsion-free}. 

\begin{lemm}
Let $\mathfrak D$ be a division algebra of degree $5$ with center 
$\ell = \bQ({\sqrt{-a}})$, where $a$ is a square-free positive integer different from $11$. Let $\tau$ be an involution of $\mathfrak D$ of 
the second kind. Then the 
subgroup $H$ of ${\mathfrak D}^{\times}$ consisting of the elements $x$ 
such that 
$x\tau(x) =1$, and ${\rm Nrd}\,(x)=1$, is torsion-free. 
\end{lemm}

\ni{\it Proof.} Let $x\in H$ be a nontrivial element of finite order. Since the 
reduced norm of $-1$ in $\mathfrak D$ is $-1$, $x\ne -1$, and therefore the 
$\bQ$-subalgebra $K := \bQ[x]$ of $\mathfrak D$ generated by $x$ is a 
nontrivial 
field extension of $\bQ$. If $K =\ell$, then $x$ lies in the center of 
$\mathfrak D$, and hence it is of order $5$. However, a nontrivial 
fifth-root of 
unity cannot be contained in a quadratic extension of $\bQ$ and so we 
conclude that $K\ne \ell$. Then $K$ is an extension of $\bQ$ of degree $5$ 
or $10$. As no extension of $\bQ$ of degree $5$ contains a root of unity 
other than $-1$, $K$ must be of degree $10$, and hence, in particular, it 
contains $\ell = \bQ({\sqrt{-a}})$. Now we note that the only roots of unity 
which can be contained in an extension of $\bQ$ of degree $10$ are the 
$11$th and $22$nd roots of unity. But the only quadratic extension contained 
in the field extension generated by either of these roots of unity is 
$\bQ(\sqrt{-11})$. Since $K\supset \bQ(\sqrt{-a})$, 
and, by hypothesis, $a\ne 11$, we have arrived at a contradiction.
\vskip1mm

\ni{\bf 5.2.} We shall now compute the covolume and the Euler-Poincar\'e 
characteristic of the principal
arithmetic subgroup $\Lambda$. 
\vskip1mm

Let $X$ be the symmetric space of $G(\bR)$, $X_u$ be the compact dual of $X$, and $\mathfrak{F}= X/\Lambda$. We note that if $m=1$, $G(\bR) = {\rm SU}(4,1)$ and $X_u =\bP_{\bC}^4$; if $m =2$, $G(\bR) = {\rm SU}(3,2)$ and $X_u ={\bf Gr}_{2,5}$; $\mathfrak F$ is a connected smooth complex projective variety.
\vskip1mm

The volume formula (5) for $n=5$ and $k =\bQ$, with the value of the Euler-product $\cE$ 
determined in 2.7, gives us 
\begin{eqnarray}
\mu(G(\bR)/\Lambda)&=& D_{\ell}^{7}\cdot\prod_{j=1}^{4}\frac{j!}{(2\pi)^{j+1}}\cdot\cE\nonumber\\
&=& D_\ell^{7}\cdot\prod_{j=1}^{4}\frac{j!}{(2\pi)^{j+1}}
\cdot \zeta(2)L_{\ell|\bQ}(3)\zeta (4)L_{\ell|\bQ}(5)\cdot\prod_{v\in\cT}e^{\prime}(P_v)\nonumber.
\end{eqnarray}

\ni From the functional equation for the $L$-function we obtain 
\begin{eqnarray*}
L_{\ell|\bQ}(3)&=&-2\pi^3 D_{\ell}^{-5/2}L_{\ell|\bQ}(-2),\\
L_{\ell|\bQ}(5)&=&\frac{2\pi^5}3 D_{\ell}^{-9/2}L_{\ell|\bQ}(-4).
\end{eqnarray*}
The following values of zeta and $L$-functions have been obtained using the software PARI.
$$\zeta(2)=\frac{\pi^2}6,\ \ \zeta(4)=\frac{\pi^4}{90},\ \ 
L_{\ell|\bQ}(-2)=-\frac{16}7,\ \ L_{\ell|\bQ}(-4)=32.$$ 
Note also that for the subgroup $\Lambda$ under consideration, $\cT = \{ 2\}$,
and $d_2 = 5$, so that (2.3) 
$$\prod_{v\in\cT}e^{\prime}(P_v)= \frac{\prod_{j=1}^{5}(2^j
  -1)}{(2^{5}-1)}= \prod_{j=1}^{4}(2^j -1).$$ 
Substituting all this in the above, we obtain 
$$
\mu(G(\bR )/\Lambda)=\prod_{j=1}^{4}\frac{j!}{(2\pi)^{j+1}}
\cdot \big(\frac{(-4\pi^{14})}{3\cdot6\cdot90}\cdot L_{\ell|\bQ}(-2)L_{\ell|\bQ}(-4)\big)
\cdot\prod_{j=1}^{4}(2^j -1)$$ 
$$=\frac9{512\pi^{14}}\cdot \frac{(-4\pi^{14})}{3\cdot 6\cdot
  90}\cdot(-\frac{16}7)\cdot32\cdot 315 \ = 1.$$ 
\vskip1mm

\ni Therefore, $\chi(\Lambda)=\chi(X_u).$ 
Theorem 3.2 of [Cl] implies that 
$H^j(\Lambda,\bC)$ 
vanishes for all odd $j$. Also, there is a natural embedding of $H^*(X_u,\bC)$ in $H^*(\Lambda,\bC)$; see [B], 3.1 and 10.2. 
Now since $\chi({\mathfrak F}) = \chi(\Lambda) =\chi(X_u)$, and for all odd $j$, $H^j({\mathfrak F}, \bC)$ (= $H^j(\Lambda, \bC)$) vanishes, 
we conclude that $\mathfrak F$ is an {\it  arithmetic fake $\bP_{\bC}^4$} if  $m=1$, and is an {\it arithmetic fake ${\bf Gr}_{2,5}$} if  $m =2$. Thus we have proved the following.

\begin{theo}
There are at least four arithmetic fake $\bP_{\bC}^4$, and at least four arithmetic fake ${\bf Gr}_{2,5}$. There does not exist any arithmetic 
fake projective space of dimension $>4$, or an arithmetic fake ${\bf Gr}_{m,n}$, with $n>5$ odd.   
\end{theo}
\vskip1mm

We next prove the following interesting result.
\begin{theo} The first integral homology group of any arithmetic fake $\bP_{\bC}^4$, and any arithmetic fake ${\bf Gr}_{2,5}$, is nonzero.\end{theo}

\ni{\it Proof.} Let $\mathfrak F$ be either an arithmetic fake $\bP_{\bC}^4$, or an arithmetic fake ${\bf Gr}_{2,5}$. Let $\Pi$ be its fundamental group. Then $H_1({\mathfrak F}, \bZ) = \Pi /[\Pi,\Pi]$. 
\vskip1mm

It follows from Theorem 2 that $\Pi$ is a cocompact torsion-free arithmetic subgroup of ${\overline G}(\bR)$, where $G$ is as in 5.1, with $m=1$ if $\mathfrak F$ is an arithmetic fake $\bP^4_{\bC}$, and $m=2$ if $\mathfrak F$ is an arithmetic fake ${\bf Gr}_{2,5}$, and $\overline G$ is the adjoint group of $G$. Proposition 1.2 of [BP] implies that $\Pi$ is actually contained in ${\overline G}(\bQ)$. We will view it as a subgroup of ${\overline G}(\bQ_2)$.   
\vskip1mm

Let $\cD$ and $\sigma$ be as in 5.1. Since ${\bQ}_2\otimes_{\bQ}\cD =({\bQ}_2\otimes_{\bQ}\ell )\otimes_{\ell}\cD={\mathfrak D}\oplus {\mathfrak D}^o$, where ${\mathfrak D}$ is a division algebra with center $\bQ_2$, 
of degree $5$, ${\mathfrak D}^o$ is its opposite, and $\sigma({\mathfrak D}) = {\mathfrak D}^o$, $G(\bQ_2)$ equals the group ${\rm SL}_1({\mathfrak D})$ of elements of reduced norm $1$ in $\mathfrak D$,  and 
${\overline G}(\bQ_2)$ equals ${\mathfrak D}^{\times}/\bQ^{\times}_2$.  We now observe that ${\overline G}(\bQ_2) = {\mathfrak D}^{\times}/ \bQ^{\times}_2$ is a pro-solvable group, i.e., if we define the decreasing sequence $\{ \cG_i\}$ of subgroups of $\cG := {\overline G}(\bQ_2)$ inductively as follows: $\cG_0 =\cG$, and $\cG_i = [\cG_{i-1},\cG_{i-1}]$, then $\bigcap \cG_i$ is trivial; to see this use Theorem 7(i) of [Ri]. This implies that for any subgroup $\cH$ of $\cG$, $[\cH,\cH ]$ is a proper subgroup of $\cH$. In particular $[\Pi,\Pi]$ is a proper subgroup of $\Pi$. This proves the theorem. 
\vskip7mm

\begin{center}
\ni{\bf 6. Five irreducible arithmetic fake $\bP^2_{\bC}\times \bP^2_{\bC}$}
\end{center}
\vskip4mm

We will now use certain results and computations of [PY] to construct five irreducible arithmetic fake $\bP^2_{\bC}\times \bP^2_{\bC}$. Let $\zeta_3$ be a 
primitive cube-root of unity, and let the pair $(k,\ell)$ of number fields be one of the following three:
$$\mathscr{C}_2 = (\bQ(\sqrt{5}), \bQ(\sqrt{5},\zeta_3)),$$
$$\hskip1cm \mathscr{C}_{10} =(\bQ(\sqrt{2}), \bQ(\sqrt{-7+4\sqrt{2}}\:)),$$
$$\mathscr{C}_{18}= (\bQ(\sqrt{6}), \bQ(\sqrt{6},\zeta_3)).$$
Let $\mathfrak v$ be the unique place of $k$ lying over $2$ if the pair is $\mathscr{C}_2$ or $\mathscr{C}_{10}$, and the unique place 
of $k$ lying over $3$ if the pair is $\mathscr{C}_{18}$. For a given pair $(k,\ell)$, let $q_{\mathfrak v}$ be the cardinality of the residue field 
of the completion $k_{\mathfrak v}$ of $k$ at $\mathfrak v$.
\vskip1mm

Let $\cD$ be a cubic division algebra with center $\ell$ whose local invariants at the two places of $\ell$ lying over $\mathfrak v$ are nonzero and 
negative of each other, and whose local invariants at all the other places of $\ell$ is zero. Then 
$k_{\mathfrak v}\otimes_{k}\cD =(k_{\mathfrak v}\otimes_k\ell)\otimes_{\ell}\cD ={\mathfrak D}\oplus {\mathfrak D}^o$, 
where $\mathfrak D$ is a cubic division algebra with center $k_{\mathfrak v}$, and ${\mathfrak D}^o$ is its opposite.   
$\cD$ admits an involution of the second kind with $k$ being the fixed field in $\ell$. We fix an involution $\sigma$ of $\cD/k$ of the second kind so 
that if $G$ is the simple simply connected $k$-group with $$G(k) = \{ z\in \cD^{\times}\:|\: z\sigma(z)=1\ \  \text{and} \ \ \mathrm{Nrd}(z)=1\},$$ 
then $G(k_v)\cong \mathrm{SU}(2,1)$ for every real place $v$ of $k$. 
\vskip1mm

As $\sigma(\mathfrak{D})= \mathfrak{D}^o$, $G(k_{\mathfrak v})$ is the compact group $\mathrm{SL}_1(\mathfrak {D})$  of elements of reduced 
norm $1$ in $\mathfrak{D}$. Let $(P_v)_{v\in V_{f}}$, be a coherent collection of maximal parahoric subgroups $P_v$ of $G(k_v)$, $v\in V_f$, 
such that $P_v$ is hyperspecial whenever $G(k_v)$ contains such a subgroup. Let $\Lambda = G(k)\cap\prod_{v\in V_f}P_v$. Let $v'$ and $v''$ 
be the two real places of $k$ and let $\mathcal{G}=G(k_{v'})\times G(k_{v''})$. Then $\mathcal{G}\cong \mathrm{SU}(2,1)\times \mathrm{SU}(2,1)$. 
Let $\overline{\mathcal{G}}$ be the adjoint group of $\mathcal{G}$. Let $X$ be the symmetric space of $\mathcal{G}$ and $X_u$ its compact dual. 
Then $X_u = \bP^2_{\bC}\times \bP^2_{\bC}$, and hence, $\chi(X_u) = 9$.  
\vskip1mm

We will view $\Lambda$ as a diagonally embedded arithmetic subgroup of 
$\mathcal{G}$. Then, in terms of the normalized Haar measure $\mu$ on 
$\mathcal{G}$ used in [P], we see using the volume formula given in that 
paper (see, also, [PY], 2.4, 2.11) that $\mu(\mathcal{G}/\Lambda) = \mu e'(P_{\mathfrak v})= 
\mu (q_{\mathfrak v}-1)^2(q_{\mathfrak v}+1)$, where the values of $\mu$ and 
$q_{\mathfrak v}$ are as given in the table in section 9.1 of [PY]. 
Moreover, according to the result in section 4.2 of [BP], the orbifold 
Euler-Poincar\'e characteristic $\chi(\Lambda)$ of $\Lambda$ equals 
$\chi(X_u)\mu(\mathcal{G}/\Lambda)= 
9\mu (q_{\mathfrak v}-1)^2(q_{\mathfrak v}+1)$. Now using the values of 
$\mu$ and $q_{\mathfrak v}$ given in the table in section 9.1 of [PY] we 
find that $\chi(\Lambda)= 3$ if $(k,\ell)$ is either $\cC_{2}$ or $\cC_{18}$, 
and $\chi(\Lambda) = 9$ if $(k,\ell) = \cC_{10}$. 
\vskip1mm

We now observe that Lemma 9.2 of [PY] holds for the group $G(k)$ described above (the proof of the lemma given in [PY] remains valid), i.e., 
$G(k)$ is torsion-free if $(k,\ell) = \cC_{10}$, and in case 
$(k,\ell)$ is either $\cC_2$ or $\cC_{18}$, the only nontrivial elements of finite order of $G(k)$ are central, and hence are of order $3$. 
Let $\overline\Lambda$ be the image of $\Lambda$ in $\overline{\mathcal{G}}$. Then $\overline\Lambda$ is a torsion-free cocompact irreducible 
arithmetic subgroup of $\overline{\mathcal{G}}$. Moreover, the natural homomorphism $\Lambda \to{\overline{\Lambda}}$ is an isomorphism if 
$(k,\ell) =\cC_{10}$, and its kernel is of order $3$ if $(k,\ell)$ is either $\cC_2$ or 
$\cC_{18}$. Hence, for each of the three pairs $(k,\ell)$, $\chi({\overline\Lambda}) = 9 =\chi(X_u)$. Let $\mathfrak{P} = X/{\overline\Lambda}$. 
Then $\mathfrak{P}$ is a smooth projective variety, and $\chi(\mathfrak{P})=\chi(\overline{\Lambda}) =\chi(X_u)$.  It is known (see the remark following 
Theorem 15.3.1 in [Ro]) that $H^j({\overline\Lambda},\bC)$ vanishes for all odd $j$. Also, there is a natural embedding of $H^*(X_u,\bC)$ in 
$H^*(\overline\Lambda,\bC)$, [B], 3.1 and 10.2.  
As $\chi(\mathfrak{P}) = \chi(X_u)$, we conclude that the Betti numbers of $\mathfrak{P}$ are same as that of $X_u =\bP^2_{\bC}\times \bP^2_{\bC}$, 
and hence $\mathfrak{P}$ is an irreducible arithmetic fake $\bP^2_{\bC}\times\bP^2_{\bC}$.
\vskip2mm

\ni{\bf Remark.} Theorem 10.1 of [PY] holds for $\mathfrak P$ (with the same proof as in [PY]), i.e., $H_1(\mathfrak{P},\bZ)$ is nontrivial.
\vskip2mm

If $(k,\ell)$ is either $\cC_2$ or $\cC_{10}$, there is a unique nonarchimedean place, say  $v_o$, of $k$, which ramifies in $\ell$. In $G(k_{v_o})$, 
up to conjugacy, there are two maximal parahoric subgroups, and thus we get four distinct $\Lambda$s. On the other hand, if $(k,\ell)=\cC_{18}$, as $D_{\ell} = D_k^2$, every place of $k$ is unramified in $\ell$, and so, up to conjugacy, we get only one $\Lambda$. Thus all together we obtain five distinct $\mathfrak{P}$s from the above construction and we have proved the following:  
\begin{theo} 
There exist at least five distinct irreducible arithmetic 
fake $\bP^2_{\bC}\times \bP^2_{\bC}$.
\end{theo}

\vskip2mm

\ni{\bf Acknowledgements.} The first-named author would like to thank the Institute for 
Advanced Study (Princeton), the Tata Institute of Fundamental Research (Mumbai), and the California Institute of Technology (Pasadena) where this work was done, for their 
hospitality. The second-named author would like to thank the University of Hong Kong, the Korea Institute for Advanced Study (Seoul) 
and Osaka University for their hospitality during the preparation of this work. The first author was supported by the Humboldt Foundation and the NSF. 
The second author received partial support from the NSF. 
\vskip1mm

We are grateful to Gunter Malle who not only supplied the list of totally complex octic fields appearing 
in the proof of Theorem 1 but also checked all the assertions about totally real number fields of degree $\leqslant 4$ 
and their totally complex quadratic extensions used in that proof. We are also grateful to Francisco Diaz y Diaz 
and J\"urgen Kl\"uners for their help with number fields. We would like to thank Don Blasius, Alireza Salehi Golsefidy, Dipendra Prasad, V.\,\,Srinivas, 
Tim Steger, Shigeaki Tsuyumine, T.N.\,\,Venkataramana and J.K.\,\,Yu for helpful conversations and correspondence. We thank P.\,\,Deligne, J.\,\,Koll\'ar and 
Y-T.\,\,Siu for their interest in the problems considered in this paper.  

\vskip3mm

\centerline{\bf References}
\vskip4mm


\vskip1.5mm
\ni[B] Borel,A., {\it Stable real cohomology of arithmetic groups}, 
Ann.\,Sci.\,Ec.\,Norm.\,Sup.\,(4) {\bf 7}(1974), 235-272. 
\vskip1.5mm
 
 \ni[BdS] Borel, A., de Siebenthal, J., {\it Les sous-groupes ferm\'es connexes de rang maximum des groupes de Lie clos}, Comm.\,Math.\,Helv. {\bf 23}(1949), 200--221.
 \vskip1.5mm
  
\ni[BP] Borel, A., Prasad, G., {\it Finiteness theorems for discrete subgroups of bounded covolume in semisimple groups.} Publ.\,Math.\,IHES No.\,{\bf 69}(1989), 119--171.
\vskip1.5mm


\ni[BS] Borevich, Z.\,I., Shafarevich, I.\,R., {\it Number theory.} Academic Press, New York (1966).
\vskip1.5mm


\ni[Cl] Clozel, L., {\it On the cohomology of Kottwitz's arithmetic varieties.} Duke 
Math.\,J. {\bf 72}(1993), 757--795.
\vskip1.5mm

\ni [Co]  Cohen, H., {\it A course in computational algebraic number theory.} Graduate Texts in Mathematics, 138. Springer-Verlag, Berlin, 1993.
\vskip1.5mm


\ni[F] Friedman, E., {\it Analytic formulas for the regulator of a number field.}
Inv.\,Math.\,{\bf 98} (1989), 599--622.
\vskip1.5mm

\ni[Kl] Klingler, B., {\it Sur la rigidit\'e de certains groupes fondamentaux, l'arithm\'eticit\'e des r\'eseaux hyperboliques complexes, et les `faux plans projectifs'}. Inv.\,Math.\,{\bf 153} (2003), 105--143.
\vskip1.5mm

\ni [M] Mumford, D., An algebraic surface with $K$ ample, $K^2=9$, $p_g=q=0.$ Amer.\,J.\,Math.\,{\bf 101}(1979), 233--244.
\vskip1.5mm

\ni[N] Narkiewicz, W., {\it Elementary and analytic theory of algebraic numbers,} third edition. Springer-Verlag, New York (2000).
\vskip1.5mm

\ni[O] Odlyzko, A.\,M., {\it  Some analytic estimates of class numbers and 
discriminants.} Inv.\,Math.\,{\bf 29} (1975), 275--286.
\vskip1.5mm

\ni[P] Prasad, G., {\it Volumes of $S$-arithmetic quotients of semi-simple groups.} Publ.\,Math. IHES No.\,{\bf 69}(1989), 91--117.
\vskip1.5mm

\ni[PY] Prasad, G., Yeung, S-K., {\it Fake projective planes}. Inv.\,Math.\,{\bf 168}(2007), 321-370. Addendum (to appear).
\vskip1.5mm

\ni[Ri] Riehm, C., {\it The norm 1 group of $\mathfrak p$-adic division algebra.} Amer.\,J.\,Math.\,{\bf 92}(1970), 499-523.
\vskip1.5mm

\ni [Ro] Rogawski, J., {\it Automorphic representations of unitary groups in three variables}, Annals of 
Math.\,Studies {\bf 123}, Princeton U.\,Press, Princeton (1990).
\vskip1.5mm

\ni[Se1] Serre, J-P., {\it Cohomologie des groupes discrets,} in Annals of 
Math.\,Studies {\bf 70}. Princeton U.\,Press, Princeton (1971).
\vskip1.5mm
 
\ni[Se2] Serre, J-P., {\it A course in arithmetic.} Springer-Verlag, New York (1973). 
\vskip1.5mm
 
\ni[T1] Tits, J., {\it Classification of algebraic semisimple groups.} Algebraic Groups and Discontinuous
Subgroups. Proc.\,A.M.S.\,Symp.\,Pure Math. {\bf 9}(1966) pp. 33--62. 
\vskip1.5mm

\ni[T2] Tits, J., {\it Reductive groups over local fields.} Proc.\,A.M.S.\,Symp.\,Pure 
Math. {\bf 33}(1979), Part I, 29--69.
 \vskip1.5mm

\ni[Y] Yeung, S.-K., {\it Integrality and arithmeticity of co-compact lattices corresponding to certain complex two-ball quotients of Picard number one}. Asian J.\,Math.\,{\bf 8}(2004), 107--130. 
\vskip1.5mm

\ni[Z] Zimmert, R., {\it Ideale kleiner Norm in Idealklassen und eine 
Regulatorabsch\"atzung.} Inv.\,Math. {\bf 62}(1981), 367--380.
\vskip1.5mm

\ni[1] The Bordeaux Database, Tables t40.001, t44.001, t60.001 obtainable from:

    ftp://megrez.math.u-bordeaux.fr/pub/numberfields/.
\vskip1.5mm



\vskip.5cm

\ni{\sc University of Michigan, Ann Arbor, MI 48109}

\ni e-mail: gprasad@umich.edu

\vskip.5cm
\ni{\sc Purdue University, West Lafayette, IN 47907}

\ni{email: yeung@math.purdue.edu}
\end{document}